\documentclass[twocolumn]{autart}    
\usepackage{cite} 
\usepackage{color}
\usepackage{graphicx}          
\usepackage{amsmath}
\usepackage{subfigure}
\usepackage{amssymb}
\usepackage{amsfonts}
\usepackage{eucal}
\usepackage[titletoc,title]{appendix}
\usepackage{mathptmx} 
\usepackage{epsfig}
\usepackage{epstopdf}

\usepackage{tikz} \usetikzlibrary{calc} \tikzset{>=latex}

\begin{document}
        
        \newtheorem{remm}{Remark}       
        \setlength\abovedisplayskip{7pt}
        \setlength\belowdisplayskip{7pt}
        \setlength\abovedisplayshortskip{7pt}
        \setlength\belowdisplayshortskip{7pt}
        \allowdisplaybreaks
        \setlength{\parindent}{1em}
        \setlength{\parskip}{0em}  
        \addtolength{\oddsidemargin}{3pt}      
        
        \begin{frontmatter}
                
                \title{ Delay-Adaptive  Compensator for 3-D Space Formation of Multi-Agent Systems with Leaders Actuation} 
                
                \thanks[footnoteinfo]{Corresponding author: Jie~Qi. The work was supported by the National Natural Science Foundation of China (62173084, 62303324, 61773112). Mamadou Diagne appreciates the support of the National Science Foundation of the USA under the grant number (FAIN): 2222250.}
                
                \author[Wang]{Shanshan Wang}\ead{wss\_dhu@126.com},  
                \author[UMICH]{Mamadou Diagne}\ead{mdiagne@ucsd.edu},
                \author[Qi]{Jie Qi}\ead{jieqi@dhu.edu.cn} 
                
                \address[Wang]{Department of Control Science and Engineering, University of Shanghai for Science and Technology, Shanghai, P. R. China, 200093}  
                \address[Qi]{College of Information science and Technology, Donghua University, Shanghai,  P. R. China, 201620} 
                \address[UMICH]{Department of Mechanical and Aerospace Engineering, University of California San Diego, La Jolla, CA, USA, 92093} 
                
                

                \begin{abstract}                          
                       \textcolor{black} {This paper focuses on the control of collective dynamics in large-scale multi-agent systems (MAS) operating in a 3-D space, with a specific emphasis on compensating for the influence of an unknown delay affecting the actuated leaders.} The communication graph of the agents is defined on a mesh-grid 2-D cylindrical surface. We model the agents' collective dynamics by a \textcolor{black}{complex-} and a real-valued reaction-advection-diffusion 2-D partial differential equations (PDEs) whose states represent the 3-D position coordinates of the agents.
                        The leader agents on the boundary suffer unknown actuator delay due to the \textcolor{black}{cumulative computation and information transmission time.} 
                        We design a delay-adaptive controller for the 2-D PDE by using PDE \textcolor{black}{backstepping combined with a Lyapunov functional method, where the latter is employed to design an  update law that generates real-time estimates of the unknown delay.}
Capitalizing on our recent result on the control of 1-D parabolic PDEs with unknown input delay, we use Fourier series expansion to bridge the control of 1-D PDEs to that of 2-D PDEs. To design the update law for the 2-D system, a new target system is defined to establish the closed-loop local  boundedness of the system trajectories in $H^2$ norm 
                        and \textcolor{black}{the regulation of the states to zero} \textcolor{black}{assuming a measurement of the spatially distributed plant's state. We illustrate the performance of   delay-adaptive controller by numerical simulations.} 
                \end{abstract}
                
        \end{frontmatter}
        
        \section{Introduction}
        \textcolor{black}{Cooperative formation control in multi-agent systems (MAS) has garnered substantial interest due to its wide-ranging applications in various engineering domains, such as UAV formation flying \cite{Mora2015}, multi-robot collaboration \cite{alonso2019distributed, Wang2017MAS}, vehicle queues \cite{Fax2004}, and satellite clusters \cite{Zetocha2000}. In MASs, communication delays, stemming from information exchange between agents, and input delays, arising from the processing/aquisition of data to update feedback control signals, can frequently lead to ``suboptimal" performance and, in more critical cases, potentially result in  system instability.   \textcolor{black}{Over the past few decades, there has been a significant body of research in multi-agent systems  focusing on communication delays. This research has predominantly employed high-order models and consensus protocols \cite{hou2017consensus, yu2010some}. In the context of non-uniform communication delays, \cite{Lee2006} establish the critical role of a globally reachable node in the information graph when designing linear agreement protocols for agents \cite{Lee2006}.} The authors of  \cite{Tian2008}  employ frequency domain analysis to derive  a delay-dependent consensus condition for a first-order multi-agent system with input and communication delays.
        In \cite{zhu2015event}, a triggering mechanism is introduced to establish a necessary and sufficient condition for leader-following consensus in multi-agent systems with input delays, while a comparable condition for second-order consensus in multi-agent dynamical systems with input delays is introduced in \cite{yu2010some}. Using the Artstein-Kwon-Pearson reduction method to convert delay-dependent
        systems into delay-free systems,   fixed-time event-triggered consensus  for 
        linear MAS with input delay is achieved in \cite{ai2021distributed}. Based on a Lyapunov method for a mean square consensus
        problem of leader-following stochastic MAS with input time-dependent or constant delay, \cite{tan2017leader} provide
        sufficient conditions to achieving consensus. A solution for leader-follower consensus in nonlinear multi-agent systems with unknown nonuniform time-varying input delay is provided in \cite{9107091} by constructing  a delay-independent output-feedback controller for each follower.
        While the prevalent focus in the literature has been on the impact of input delay on follower agents,  \cite{qi2019control} address a known delay affecting the actuated leaders within a 3-D infinite-dimensional framework. Furthermore, most of these studies rely on ordinary differential equations (ODEs) models, namely, each agent's dynamic state is represented by an ODE, resulting in increased system complexity as the number of agents grows \cite{Lee2006,Lin2014}.}

         \textcolor{black}{For multi-agent systems, control designs using partial differential equations (PDEs) provide a compact representation for capturing the dynamics of large-scale systems. These PDEs, whether they take a parabolic or hyperbolic form, describe the position coordinates of individual agents, as demonstrated in various works including \cite{Meurer2011, FREUDENTHALER2020108897,   Frihauf2011, Qi2015, qi2017wave} and the reference therein. In the case of  parabolic systems,  the diffusion term, namely, the Laplace operator plays the role of MAS consensus protocol modeled by ODEs. Actuation of the leader agents positioned on the periphery of the communication structure demands a greater amount of information and computational resources compared to the follower agents. Consequently, leaders are more susceptible to delays that affect the formation control. Using the nominal delay-compensated boundary control law proposed in \cite{MKrstic2009,Wang2017}, the authors of \cite{qi2019control} designed  a boundary feedback law for MAS in 3-D space under a constant and known input delay. However, in practical scenarios, knowing precisely the value of the delay is often unfeasible, and instead, it is possible to estimate only its upper and lower bounds.  To overcome such a challenge, the authors of \cite{LIU2021104927} investigate the determination of the delay bounds within which regulated state synchronization is attainable for a multi-agent system with unknown and nonuniform input delays. Similarly, in \cite{Zhangunknowndelay2021}, such a delay bound is characterized for semi-global state synchronization in a multi-agent system with actuator saturation and unknown nonuniform input delays. Nevertheless, there is a dearth of literature that addresses the issue of unknown delays in the context of multi-agent systems modeled by partial differential equations (PDEs). \textcolor{black}{For a reaction-diffusion systems subject to an unknown boundary input delays \cite{WANG2021109909} pioneering exploration led to a delay-adaptive compensated controller that ensures \textcolor{black}{the regulation of the system's state to zero.} Motivated by decontamination of a polluted surface, \cite{wang2021delay} constructed a delay-adaptive predictor feedback for reaction-diffusion systems subject to a delayed distributed input. The stabilization of deep-sea construction vessels using Batch-Least Squares Identifiers \cite{karafyllis2019adaptive} has been achieved in \cite{Jiwang2023delay} where finite-time exact identification of an unknown boundary input delay and simultaneously exponential regulation of the plant's state for a hyperbolic PDE-ODE system is ensured.  More recently, a Lyapunov design approach that enables \textcolor{black}{global} stability for a hyperbolic PIDE (Partial Integro-Differential Equation) with an unknown boundary input delay was introduced in \cite{wang2023delay}.} }

        We consider a formation control in 3-D space of a multi-agent system with  unknown actuator delay. The collective dynamics of the MAS is modeled by two diffusion 2-D PDEs; one is a complex-valued PDE whose states represent the agents' positions in coordinates $(x,y)$ and the other is a real-valued PDE whose states represent the agents' positions in coordinate $z$. \textcolor{black}{We utilize PDE backstepping design in conjunction with a Lyapunov method to construct a dynamic delay-adaptive boundary feedback law.  The nominal backstepping controller acquires complementary information about the unknown parameter through an update law driven by a carefully designed ODE. Our present work potentially marks a pioneering contribution in the field of PDE-based formation control for multi-agent systems with unknown input delays in three-dimensional space.}  We introduce Fourier series expansion to reduce the dimensionality of the 2-D system to $n$ 1-D systems. In contrast to the result in \cite{qi2019control}, the   target system in the present study accounts for several highly nonlinear terms, generated by the delay-adaptive scheme, which pose challenges in establishing the convergence of their series representations, a crucial prerequisite for transforming the 1-D  system into a 2-D  system.


        
        This paper is organized as follows.  Section \ref{Multi-agent} introduces the PDE-based model for a MAS with actuation delay. Section \ref{controller} presents the delay-adaptive control design for the MAS collective dynamics subject to unknown actuation delay. \textcolor{black}{The design of the delay's adaptation law  and  the local stability analysis of the MAS are presented  in Section \ref{adaptive} and Section \ref{5-proof}, respectively. Simulation results are provided in Section \ref{simu}.  The paper concludes with a discussion possible of   future works  in  Section \ref{perspec}.}
        
        \textbf{Notation:} Throughout the paper, we adopt the following notation for  $\chi(\theta,s)$ on the cylindrical surface \cite{tang2017formation,qi2019control}:

        The $L^2$ norm is defined as      
        \begin{align}
                \label{equ-define11} 
                ||\chi(s,\theta)||_{L^2}=  
                \left(\int_{0}^{1}\int_{-\pi}^{\pi}|\chi(s,\theta)|^2
                \mathrm{d}\theta\mathrm{d}s\right)^\frac{1}{2},
        \end{align}
        for $\chi(\theta,s)\in L^{2}((0,1)\times(-\pi,\pi))$. To save space, we  set $\rVert  \chi\rVert^2=\rVert \chi(s,\theta)\rVert^2_{L^2}$.       
        The Sobolev norm $||\centerdot||_{H^1}$ is defined as  $||\chi||^2_{H^1}
                =\|\chi\|^{2}+
                \|\partial_s\chi\|^{2}
                +\|\partial_\theta\chi\|^{2}$
        for $\chi(\theta,s)\in H^{1}((0,1)\times(-\pi,\pi))$.        
        The Sobolev norm $||\centerdot||_{H^2}$ is defined as  $||\chi||^2_{H^2}=\|\chi\|_{H^1}^{2}+\|\partial^2_{s}\chi\|^{2}+2\|\partial_{s\theta}\chi\|^{2}+\|\partial^2_{\theta}\chi\|^{2}$
for $\chi(\theta,s)\in H^{2}((0,1)\times(-\pi,\pi))$.

        \section{\protect Muilt-Agent's PDEs Model}\label{Multi-agent}
        \subsection{Model description}
\begin{figure}
                \centering
                \includegraphics[width=0.4\textwidth]{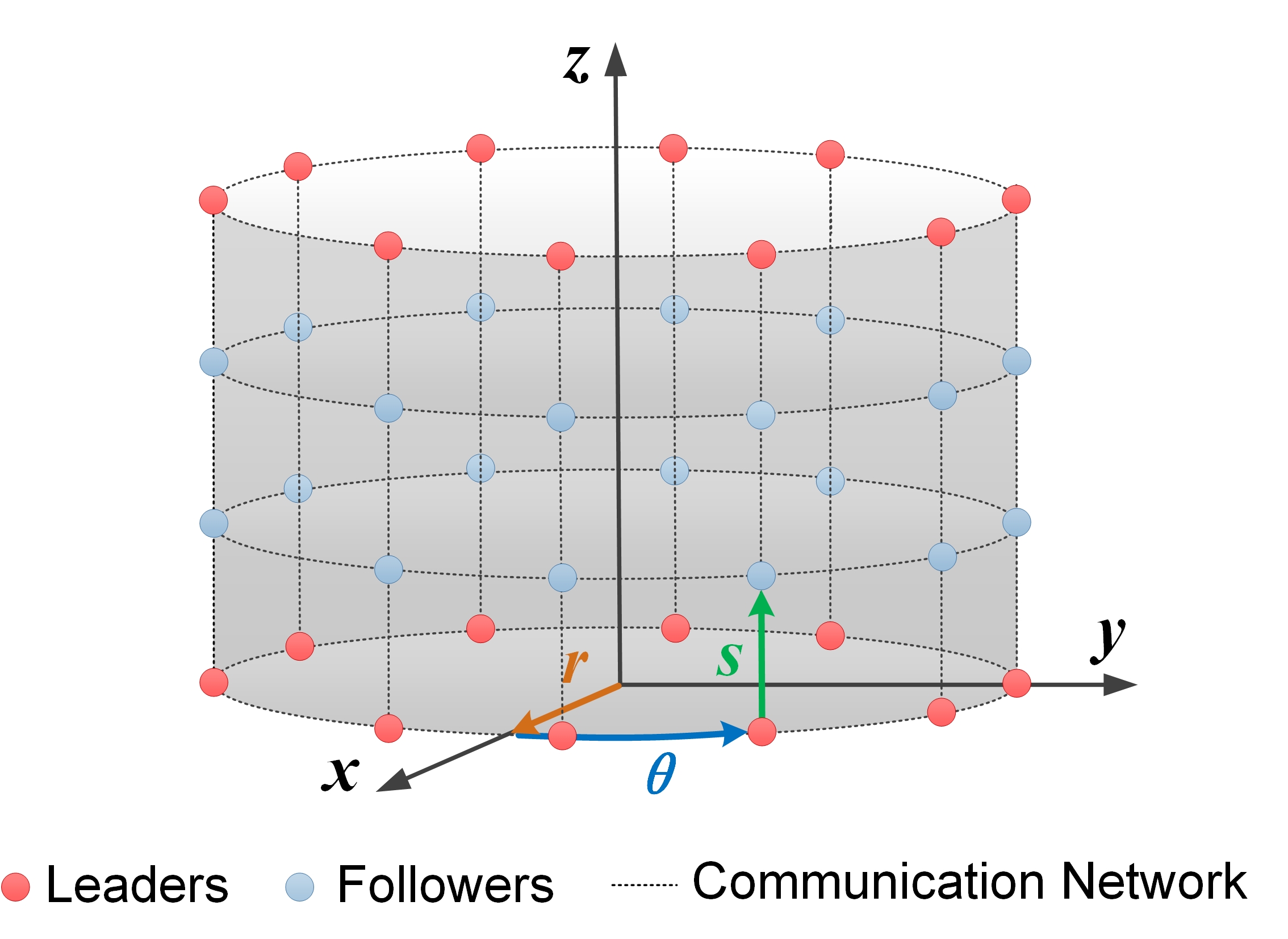}
                \vspace{-0.21cm}
                \caption{Cylindrical surface topology prescribing the communication relationship among agents.  The agents at the uppermost and lowermost layers are leaders. Each follower has four neighbors.} \label{figure-topology}
        \end{figure}

        Following \cite{qi2019control},  we consider a group of agents located on a cylindrical surface undirected topology graph with index $(i,j),~i=1,..., \ M,~j=1,..., N$, moving in a 3-D space  under the coordinate axes $(x,y,z)$.  
        A complex-valued state $u=x+\mathrm{j}y$ is defined to simplify the expression of the components on the  $(x,y)$ axes. 
        Defining  the discrete indexes $(i,j)$ of the agents into a continuous domain, 
        $\Omega=\{(s, \theta):0<{s}<{1}$, $-\pi\leq{\theta}\leq{\pi}\}$, as $M, ~N\to {\infty}$ (see, Figure \ref{figure-topology}), the continuum model of  the collective dynamics of a large scale multi-agent system as follows 
        \begin{align}
                \label{equ-u_0}
                &{\partial}_{t}u(s,{\theta},t)={\Delta}{u}(s,{\theta},t)
                +\beta_1{\partial}_su(s,{\theta},t)+{\lambda_{1}}{u}(s,{\theta},t),
                \\
                \label{equ-z_0}
                &{\partial}_{t}{z}(s,{\theta},t)={\Delta}{z}(s,{\theta},t)
                +\beta_2{\partial}_sz(s,{\theta},t)+{\lambda_{2}}{z}(s,{\theta},t),
                \\
                &\textcolor{black}{{u}(s,-\pi,t)={u}(s,\pi,t)},\quad {u}(0,{\theta},t)=f_1(\theta),\label{equ-bdu}\\
                &{u}(1,{\theta},t)=g_1(\theta)+{U}({\theta},t-D),
                \label{bnd-U}
                \\
                \label{equ-bdz}
                &\textcolor{black}{{z}(s,-\pi,t)={z}(s,\pi,t)},\quad{z}(0,{\theta},t)=f_2(\theta),\\
                \label{equ-bdz1}
                &{z}(1,{\theta},t)=g_2(\theta)+{Z}({\theta},t-D), 
        \end{align}
        where $(s,{\theta},t) \in\Omega \times {\mathbb{R}^+}$, $u$, $\lambda_{1}$, $ \beta_1 \in{\mathbb{C}}$, $z$, $\lambda_2$, $\beta_2 \in{\mathbb{R}}$. The coordinates $(s,\theta)$ are the spatial variables denoting the indexes of the agents in the continuum and ${\Delta}$ represents the following  Laplace operator
\begin{align}
\label{equ-lapu}
{\Delta}{u}(s,{\theta},t) & ={\partial}^2_{s}u(s,{\theta},t)+{\partial}^2_{\theta}u(s,{\theta},t), \\
\label{equ-lapz}
{\Delta}{z}(s,{\theta},t) & ={\partial}^2_{s}z(s,{\theta},t)+{\partial}^2_{\theta}z(s,{\theta},t),
\end{align}
        which is defined as ''consensus operators" for  PDE representations \cite{Ferrari2006}. Note that 
        the boundary conditions (\ref{equ-bdu}) and (\ref{equ-bdz}) are periodical on the cylinder surface  (see Figure \ref{figure-topology}) while $f_1(\theta)$, $g_1(\theta)$, $f_2(\theta)$ and $g_2(\theta)$ are non-zero bounded boundary conditions for the states $u$ and $z$, respectively. 
        
        To control the MAS to desired formations, we consider a configuration where the agents at the boundaries $s=0$ and $s=1$ are the leaders that drive all the agents to prescribe equilibrium.  In \eqref{bnd-U} and \eqref{equ-bdz1}, we defined the input delay  $D>0$ affecting the actuated leaders and caused by communication lags in leader-follower configurations. In practice, the exact value of the delay is hard to measure, only the bounds of the unknown delay can be estimated, so we assume:
        \newtheorem{assuption}{\bf{Assumption}}
        \begin{assuption}\label{as1}
                \rm{Assume  delay $D\in \{D\in \mathbb{R}^+|\underline{D}\leq D\leq \overline {D}\}$, where $\underline{D}$ and $\overline {D}$ are the known lower and upper bounds, respectively.}
        \end{assuption}
        
        \begin{remm}\rm{Letting $\partial_t u(s,\theta,t)=0$ and $\partial_t z(s,\theta,t)=0$, one can solve \eqref{equ-u_0}--\eqref{equ-bdz1} without control and get  steady state profiles 
                        $\bar{u}(s,\theta)$ and $\bar{z}(s,\theta)$, which express the desired formations 
                        as  functions  of the values of the  parameters $\lambda_{1}$, $\beta_1$, $\lambda_2$, $\beta_2$  and the  open-loop boundary conditions   $f_1(\theta)$, $g_1(\theta)$, $f_2(\theta)$ and $g_2(\theta)$ \cite{qi2019control}. }
        \end{remm}
        
        \section{{Delay-}adaptive boundary controller design}\label{controller}
        
        First, define the  error between the actual system and the desired system as  $\tilde{u}(s,{\theta},t)=u(s,{\theta},t)-\bar{u}(s,{\theta})$,
        and then introduce  a change of variable  $\phi(s,{\theta},t)=e^{\frac{1}{2}\beta_1s}\tilde{u}(s,{\theta},t)$ for removing the convection term,
        \begin{align}
                \label{equ-phi-ori}
                &\partial_{t}{\phi}(s,{\theta},t)={\Delta}{\phi}(s,{\theta},t)+{\lambda'_1}{\phi}(s,{\theta},t),~~s\in(0,1),\\
                &\textcolor{black}{{\phi}(s,-\pi,t)={\phi}(s,\pi,t)},\quad {\phi}(0,{\theta},t)=0,\\
                &{\phi}(1,{\theta},t)=\Phi(\theta,t-D),\label{initial-phi-ori} 
        \end{align}
        where $\lambda'_1=\lambda_1-\frac{1}{4}\beta_1^2$ and ${\Phi}({\theta},t-D)=e^{\frac{1}{2}\beta_1}{U}({\theta},t-D)$.
        By employing a transport PDE of $\vartheta$ representation of the delay appearing in \eqref{bnd-U}, we transform the error system  \eqref{equ-phi-ori}--\eqref{initial-phi-ori} as follows:
        \begin{align}
                \label{equ-phi}
                &\partial_{t}{\phi}(s,{\theta},t)={\Delta}{\phi}(s,{\theta},t)+{\lambda'_1}{\phi}(s,{\theta},t),~~s\in(0,1),\\
                &\textcolor{black}{{\phi}(s,-\pi,t)={\phi}(s,{\pi},t)},\quad{\phi}(0,{\theta},t)=0, \\
                &{\phi}(1,{\theta},t)={\vartheta}(0,{\theta},t),\\
                &D\partial_{t}{\vartheta}(s,{\theta},t)=\partial_{s}{\vartheta}(s,{\theta},t),\label{equ-vartheta-0}\\
                &\textcolor{black}{\vartheta(s,{-\pi},t)=\vartheta(s,\pi,t)},\quad\vartheta(1,{\theta},t)=\Phi(\theta,t),\label{initial-vartheta}
        \end{align}
        where $\vartheta(s,{\theta},t)=\Phi({\theta},t+D(s-1))$, defined in $\Omega\times\mathbb{R}^+$.
        In the following, we will adopt the design method presented in \cite{qi2019control} to derive the dynamic boundary adaptive controller  of the states  $u$ and $z$ but limit our analysis to the $u$ component of the state as a similar approach applies for the $z$ dynamics. 
        
        
        \subsection{Fourier series expansions}
        
        In order to transform  the 2-D system 
        \eqref{equ-phi}--\eqref{initial-vartheta}   into $n$ 1-D systems,   we introduce the   Fourier series expansion  \cite{vazquez2016explicit} as
        \begin{align}
                \label{fourier-phi}
                &{\phi}(s,{\theta},t)=\sum_{n=-\infty}^{\infty}{\phi}_{n}(s,t)e^{\mathrm{j}n\theta},~~ \\
&{\Phi}({\theta},t)=\sum_{n=-\infty}
                ^{\infty}{\Phi}_{n}(t)e^{\mathrm{j}n\theta},\\
                \label{fourier-Phi}
                &{\vartheta}(s,{\theta},t)=\sum_{n=-\infty}^{\infty}{\vartheta}_{n}(s,t)e^{\mathrm{j}n\theta},
        \end{align} 
        where $\phi_n$, $\Phi_n$, ${\vartheta}_{n}$ are  the $n$ Fourier coefficients; independent of the angular argument $\theta$.  As an illustration,  one of the coefficients in \eqref{fourier-phi}--\eqref{fourier-Phi} is given as $\phi_n(s,t)=\frac{1}{2\pi}\int_{-\pi}^{\pi}\phi(s,\psi,t)e^{-\mathrm jn\psi}\mathrm d\psi$.
        Introducing \eqref{fourier-phi}--\eqref{fourier-Phi} to
        \eqref{equ-phi}--\eqref{initial-vartheta},
        we  get the following 1-D PDE of the Fourier coefficients  ${\phi}_{n}(s,t)$ and ${\vartheta}_{n}(s,t)$
        \begin{align}
                \label{phi-n}
                &\partial_t{\phi}_{n}(s,t)=\partial^2_s{\phi}_{n}(s,t)+(\lambda'_1-n^2){\phi}_{n}(s,t),~s\in(0,1),\\
                &{\phi}_{n}(0,t)=0, \quad {\phi}_{n}(1,t)={\vartheta}_{n}(0,t), \label{bnd-phi1} \\
                &D\partial_t{\vartheta}_{n}(s,t)={\partial_s}{\vartheta}_{n}(s,t), \quad {\vartheta}_{n}(1,t)={\Phi}_{n}(t).\label{initial-vartheta-n}
        \end{align} 
        We   design a feedback adaptive  controller $\Phi_n$ to stabilize each cascade system $(\phi_n, \vartheta_n)$   in \eqref{phi-n}--\eqref{initial-vartheta-n} by postulating the following  transformations 
        {\begin{align}
                        \label{equ-transformation}
                        {w}_n(s,t)&=\mathcal{T}_n[\phi_n](s,t)={\phi}_n(s,t)-\int_{0}^s
                        k_n(s,\tau)\phi_n(\tau,t)\mathrm{d}\tau,\\
                        \label{equ-transformation2}
                        \nonumber{h}_n(s,t)&=\mathcal{T}_n[\vartheta_n](s,t)={-\hat D(t)}\int_{0}^s p_n(s,\tau,\hat D(t))\vartheta_n(\tau,t)\mathrm{d}\tau\\
                        &+{\vartheta}_n(s,t)-\int_{0}^1 \gamma_n(s,\tau,\hat D(t))\phi_n(\tau,t)\mathrm{d}\tau,
                \end{align}
                with the inverse  transformations 
                \begin{align}
                        \label{equ-transformation-inverse}
                        {\phi}_n(s,t)&=\mathcal{T}_n^{-1}[w_n](s,t)={w}_n(s,t)+\int_{0}^s
                        l_n(s,\tau)w_n(\tau,t)\mathrm{d}\tau,\\
                        \label{equ-transformation2-inverse}
                        \nonumber{\vartheta}_n(s,t)&=\mathcal{T}_n^{-1}[h_n](s,t)={\hat D(t)}\int_{0}^s q_n(s,\tau,\hat D(t))h_n(\tau,t)\mathrm{d}\tau\\
                        &+h_n(s,t)+\int_{0}^1 \eta_n(s,\tau,\hat D(t))w_n(\tau,t)\mathrm{d}\tau,
        \end{align}}
        where the kernels are defined on ${\mathcal D}=\{(s,\tau,\hat D(t)):0\leq \tau \leq s \leq 1\}$, and  $\hat D(t)$ is the estimate of unknown input delay\footnote{For the sake of simplicity,    $\hat D(t)$ is defined as $\hat D$ in the remaining part of our developments.}. 
        
        Hence, by PDE Backstepping method,  \eqref{phi-n}--\eqref{initial-vartheta-n} maps into the following target system parameterized by $\hat D=\tilde D-D$ 
        \begin{align}
                &\partial_t{w}_{n}(s,t)=\partial^2_s{w}_n(s,t)-n^2{w}_n(s,t),\label{equ-w_n-0}\\
                &{w}_n(0,t)=0,\quad{w}_n(1,t)={h}_n(0,t),\\
                &D{\partial_th}_n(s,t)={\partial_sh}_n(s,t)-\tilde D P_{1n}(s,t)-D\dot{\hat D}P_{2n}(s,t),\\
                &{h}_n(1,t)=0,\label{equ-h_n-0}
        \end{align}
        where
        \begin{align}
                \nonumber P_{1n}&(s,t)=\int_{0}^1\bigg(-{\partial_\tau}\gamma_n(s,1,\hat D)l_n(1,\tau)+\frac{1}{\hat D}\partial_s \gamma_n(s,\tau,\hat D)\\
                \nonumber&+\frac{1}{\hat D}\int_{\tau}^1\partial_s \gamma_n(s,\tau,\hat D)l_{n}(\xi,t)\mathrm{d}\xi\bigg)w_{n}(\tau,t)\mathrm{d}\tau\\
                &-{\partial_\tau}\gamma_n(s,1,\hat D)h_{n}(0,t),\\
                \nonumber P_{2n}&(s,t)=\int_{0}^1\bigg(\int_\tau^1\partial_{{\hat D}} \gamma_n(s,\tau,\hat D)l_n(\xi,\tau)\mathrm d\xi+\partial_{{\hat D}} \gamma_n(s,\tau,\hat D)\\
                \nonumber&+\int_{0}^s (p_n(s,\xi,\hat D)+{\hat D}\partial_{{\hat D}} p_n(s,\xi,\hat D))\eta_n(\xi,\tau,\hat D)\mathrm{d}\xi\bigg)  \nonumber\\
                \nonumber&\cdot w_n(\tau,t)\mathrm{d}\tau+\int_{0}^s\bigg({\hat D}\partial_{{\hat D}} p_n(s,\tau,\hat D)+ p_n(s,\tau,\hat D)\\
                \nonumber&+\hat D\int_{\tau}^s (p_n(s,\xi,\hat D)+{\hat D}\partial_{{\hat D}} p_n(s,\xi,\hat D))q_n(\xi,\tau,\hat D)\mathrm{d}\xi \bigg) \\
                &\cdot  h_n(\tau,t)\mathrm d\tau.
        \end{align}
        The mapping \eqref{equ-transformation}, \eqref{equ-transformation2}  is well defined if  the kernel functions $k_n(s,\tau)$, $\gamma_n(s,\tau)$ and $p_n(s,\tau)$  satisfy
        \begin{align}
                \label{equ-kn}
                &\partial^2_s{k}_{n}(s,\tau)=\partial^2_\tau{k}_{n}(s,\tau)+\lambda'_1{k}_{n}(s,\tau),\\
                &{k}_{n}(s,0)=0,\quad {k}_{n}(s,s)=-\frac{\lambda'_1}{2}s,\label{bnd-kn}\\
                &\partial_s{\gamma}_{n}(s,\tau,\hat D)=\hat D(\partial^2_\tau{\gamma}_{n}(s,\tau,\hat D)+(\lambda'_1-n^2){\gamma}_{n}(s,\tau,\hat D)),\label{gamman}\\
                &{\gamma}_{n}(s,0,\hat D)=
                {\gamma}_{n}(s,1,\hat D)=0,  \quad{\gamma}_{n}(0,\tau,\hat D)={k}_{n}(1,\tau),\label{initial-gamman} \\
                &\partial_sp_{n}(s,\tau,\hat D)=-\partial_\tau p_{n}(s,\tau,\hat D),\label{pn}\\ &p_{n}(s,1,\hat D)=-\partial_\tau{\gamma}_{n}(s,\tau,\hat D)|_{\tau=1}.\label{bnd-pn}
        \end{align} 
        The 
        solution of the above gain kernels PDEs is given by
        \begin{align}
                \label{equ-kn1}
                &{k}_{n}(s,\tau)=-\lambda\tau\frac{I_1(\sqrt{\lambda'_1(s^2-{\tau}^2)})}{\sqrt{\lambda'_1(s^2-{\tau}^2)}},\\
                \label{equ-gamman5}
                \nonumber&\gamma_n(s,\tau,\hat D)=2\sum_{i=1}^{\infty}e^{\hat D(\lambda'_1-n^2-i^2{\pi^2})s}{\mathrm{sin}}(i\pi\tau) \int_{0}^{1}{\mathrm{sin}}(i\pi\xi)\\
                &~~~~~~~~~~~~~~~~~\cdot k(1,\xi){\mathrm{d}{\xi}},\\
                &{p}_{n}(s,\tau,\hat D)=-{\partial}_2{\gamma}_{n}(s-\tau,1,\hat D), 
        \end{align}
        where $\partial_2\gamma_{n}(\cdot,\cdot,\cdot)$  denotes the derivative of $\gamma_n(\cdot,\cdot,\cdot)$ with respect to the second argument.  {Similarly, one can get the kernels in  inverse transformations
                \eqref{equ-transformation-inverse}, \eqref{equ-transformation2-inverse}:}
        \begin{align}
                \label{equ-ln1}
                &{l}_{n}(s,\tau)=-\lambda\tau\frac{J_1(\sqrt{\lambda'_1(s^2-{\tau}^2)})}{\sqrt{\lambda'_1(s^2-{\tau}^2)}},\\
                \label{equ-eta5}
                \nonumber&\eta_n(s,\tau,\hat D)=2\sum_{i=1}^{\infty}e^{-\hat D(n^2+i^2\pi^2)s}{\mathrm{sin}}(i\pi\tau) \int_{0}^{1}{\mathrm{sin}}(i\pi\xi)\\
                &~~~~~~~~~~~~~~~~~\cdot k(1,\xi){\mathrm{d}{\xi}},\\
                &{q}_{n}(s,\tau,\hat D)=-{\partial}_2{\eta}_{n}(s-\tau,1,\hat D). 
        \end{align}
        
        From \eqref{initial-vartheta-n}, \eqref{equ-transformation2} and \eqref{equ-h_n-0}, the 1-D delay-compensated  adaptive controller writes
        \begin{align}
                \label{equ-controller}
                \nonumber{\textcolor{black}{\Phi}}_n(t)=&{\hat D}\int_{0}^1 p_n(1,\tau,\hat D)\vartheta_n(\tau,t)\mathrm{d}\tau \\
&+\int_{0}^1 \gamma_n(1,\tau,\hat D)\phi_n(\tau,t)\mathrm{d}\tau.
        \end{align}
        
        \subsection{2-D delay-compensated  adaptive controller}
        
        In order to obtain the 2-D delay-compensated adaptive  controller, we assemble all the $n$ 1-D transformations defined in \eqref{equ-transformation}--\eqref{equ-transformation2} in the form of Fourier series to recover the  2-D  domain components and then get 
        \begin{align}\label{transform-phi-w}
                \nonumber w(s,\theta,t)&{=\sum_{n=-\infty}^{\infty}{w}_{n}(s,t)e^{\mathrm{j}n\theta}}\\ &={\phi}(s,\theta,t)
                -\int_{0}^{s}{k}(s,\tau){\phi}
                (\tau,\theta,t){\mathrm{d}{\tau}},\\
                \nonumber{h}(s,{\theta},t) 
                &{=\sum_{n=-\infty}^{\infty}h_{n}(s,t)e^{\mathrm{j}n\theta}}={\vartheta}(s,{\theta},t)\\
                \nonumber&-\int_{0}^1\int_{-\pi}^\pi \gamma(s,\tau,\theta,\psi,\hat D)
                w(\tau,\psi,t)\mathrm{d}\psi\mathrm{d}\tau  \nonumber \\
                &-\hat D\int_{0}^s\int_{-\pi}^\pi p(s,\tau,\theta,\psi,\hat D)\vartheta(\tau,\psi,t)\mathrm{d}\psi\mathrm{d}\tau, \label{transform-var-h}
        \end{align}
        where $k(s,\tau)$ is defined in \eqref{equ-kn1}, the related 2-D kernels are given as 
        \begin{align}
                \label{kernel-gamma}
                \gamma(s,\tau, \theta,\psi,\hat D)=&2{Q(s,\theta-\psi,\hat D)}\sum_{i=1}^{\infty}e^{\hat D(\lambda'_1-i^2\pi^2)s}
                \mathrm{sin}(i\pi\tau)\nonumber\\
                &\cdot\int_{0}^{1}\mathrm{sin}(i\pi\xi)k(1,\xi)\mathrm{d}\xi,\\
                \label{kernel-p}
                p(s,\tau,\theta,\psi,\hat D)=&-\partial_2\gamma(s-\tau,1,\theta-\psi,\hat D).
        \end{align}
        For all $s\in[0,1] $, defining ${Q(s,\theta-\psi,\hat D)}=\frac{1}{2\pi}\sum_{n=-\infty}^{\infty}
                e^{-\hat Dn^2s}$\\$\cdot e^{\mathrm{j}n(\theta-\psi)}$.
        {Due to $0<Q(s,\theta-\psi,\hat D)\leq P(e^{-\hat Ds},\theta -\psi)$, where $P$ denotes Possion Kernel.} Using the properties of  Poisson kernels \cite{Brown2009}, one  gets  the boundedness of  the kernel functions $\gamma(s,\tau, \theta,\psi,\hat D)$ and $p(s,\tau, \theta,\psi,\hat D)$. In a similar way, we get
        the inverse  transformations of \eqref{transform-phi-w} and \eqref{transform-var-h}  are  given by 
        \begin{align}
                \label{equ-trans3}
                &{\phi}(s,{\theta},t)={w}(s,{\theta},t)+\int_{0}^s l(s,\tau)w(\tau,{\theta},t)\mathrm{d}\tau,\\
                &\nonumber{\vartheta}(s,{\theta},t)=\int_{0}^1\int_{-\pi}^\pi \eta(s,\tau,\theta,\psi,\hat D)w(\tau,\psi,t)\mathrm{d}\psi\mathrm{d}\tau\\
                \label{equ-trans4}
                &{+h}(s,{\theta},t)+\hat D\int_{0}^s\int_{-\pi}^\pi q(s,\tau,\theta,\psi,\hat D)h(\tau,\psi,t)\mathrm{d}\psi\mathrm{d}\tau,
        \end{align} 
        where the gain  kernels $l$, $q$ and $\eta$ are defined as 
        \begin{align}
                \label{equ-kn2}
                &l(s,\tau)=-\lambda\tau\frac{J_1(\sqrt{\lambda'_1(s^2-{\tau}^2)})}{\sqrt{\lambda'_1(s^2-{\tau}^2)}},\quad \\
                &\eta(s,\tau,\theta,\psi,\hat D)=2{Q(s,\theta-\psi,\hat D)}\sum_{i=1}^{\infty}e^{-\hat Di^2\pi^2s}\mathrm{sin}(i\pi\tau)\nonumber \\
                &~~~~~~~~~~~~~~~~~~~~~~~~\cdot\int_{0}^{1}\mathrm{sin}(i\pi\xi)l(1,\xi)\mathrm{d}\xi,\label{equ-rn2}\\
                &q(s,\tau,\theta,\psi,\hat D)=-{\partial_2}{\eta}(s-\tau,1,\theta,\psi,\hat D)\label{equ-qn2}.
        \end{align}

        From \eqref{equ-controller}, \eqref{kernel-gamma} and \eqref{kernel-p}, we obtain the  following delay-adaptive control law: 
        \begin{align}
                \label{equ-endU}
                &\nonumber U({\theta},t)=\int_0^1 \int_{-\pi}^\pi \gamma(1,\tau,\theta,\psi,\hat D)e^{-\frac{1}{2}\beta_1(1-\tau)}(u(\tau,{\psi},t)\\
                \nonumber&~~~~-\bar{u}(\tau,{\psi}))\mathrm{d}\psi\mathrm{d}\tau-\int_{t-D}^{t}\int_{-\pi}^\pi {\partial_2}{\gamma}
                (\frac{t-\nu}{\hat D},1,\theta,\psi,\hat D)e^{\frac{1}{2}\beta_1}\\
                &~~~~\cdot U({\psi},\nu)\mathrm{d}\psi\mathrm{d}\nu.
        \end{align}

        \subsection{Target system for the plant with  unknown input delay}
        Similarly, in order to obtain the 2-D Target system for the plant with  unknown input delay, we assemble all the $n$ 1-D target systems defined in \eqref{equ-w_n-0}--\eqref{equ-h_n-0} in the form of Fourier series to return back to the  2-D  domain 
        \begin{align}
                \label{equ-aw0-adp}
                &\partial_t{w}(s, {\theta}, t)={\Delta}{w}(s, {\theta}, t), \\
                &\textcolor{black}{{w}(s, -\pi, t)={w}(s, {\pi}, t),}  \\ 
                &{w}(0, {\theta}, t)=0, \quad {w}(1, {\theta}, t)={h}(0, {\theta}, t)\label{bnd-aw-adp}, \\
                & D\partial_t{h}(s, {\theta}, t)=\partial_s{h}(s, {\theta}, t)-\tilde DP_1(s, \theta, t)-D\dot{\hat D}P_2(s, \theta, t), 
                \label{equ-ah-adp}\\
                &\textcolor{black}{{h}(s, -\pi, t)={h}(s, {\pi}, t),}\quad {h}(1, {\theta}, t)=0,\label{initial-ah-adp}
        \end{align}
        with
        \begin{align}
                \nonumber P_{1}(s, \theta, t)=&\int_0^1\int_{-\pi}^\pi M_1(s, \tau, \hat \theta, \psi,t)w(\tau, \psi, t)\mathrm  d\psi\mathrm d\tau\\
                &+\int_{-\pi}^\pi M_{2}(s,  \theta, \psi,t)h(0, \psi, t)\mathrm d\psi, \label{equ-P1}\\
                \nonumber P_{2}(s, \theta, t)=&\int_0^1\int_{-\pi}^{\pi}M_3(s, \tau, \theta, \psi,t)w(\tau, \psi, t)\mathrm d\psi\mathrm d\tau\\
                &+\int_0^s\int_{-\pi}^{\pi}M_4(s, \tau, \theta, \psi,t)h(\tau, \psi, t)\mathrm d\psi\tau, \label{equ-P2}
        \end{align}
        where $M_{i}, ~ i=1, 2, 3, 4$ are functions defined below:
        \begin{align}
                \label{equ-M1}
                &\nonumber M_1(s, \tau, \theta, \psi,t)=\frac{1}{\hat D}\bigg(\int_{\tau}^{1}\gamma_s(s,\xi,\theta, \psi,\hat D)l(\xi, \tau)\mathrm d\xi\\
                &~~~~+\gamma_{s}(s, \tau,\theta, \psi,\hat D)\bigg)-\textcolor{black}{{\partial_2}\gamma(s, 1, \theta, \psi,\hat D)}l(1,\tau),\\
                &M_2(s,\tau,  \theta, \psi,t)=-\textcolor{black}{{\partial_2}\gamma(s, 1, \theta, \psi,\hat D)}, \\
                &\nonumber M_{3}(s, \tau, \theta, \psi,t)=\int_\tau^1\gamma_{\hat D}(s, \xi, \theta, \psi,\hat D)l(\xi, \tau)\mathrm d\xi\\
                \nonumber&~~~~+\int_0^s\int_{-\pi}^\pi(p(s, \xi, \theta, \varphi,\hat D)+\hat Dp_{\hat D}(s, \tau, \theta,  \varphi,\hat D))\\
                &~~~~\cdot\eta(\xi, \tau, \varphi, \psi,\hat D)\mathrm d \varphi\mathrm d\xi+\gamma_{\hat D}(s, \tau, \theta, \psi,\hat D),\\
                &\nonumber M_{4}(s, \tau, \theta, \psi,t)=p(s, \tau, \theta, \psi,\hat D)+\hat Dp_{\hat D}(s, \tau, \theta, \psi,\hat D)\\
                \nonumber&~~~~+\hat D\int_\tau^s\int_{-\pi}^\pi(p(s, \xi, \theta,  \varphi,\hat D)+\hat Dp_{\hat D}(s, \xi, \theta,  \varphi,\hat D))\\
                &~~~~\cdot q(\xi, \tau,  \varphi, \psi,\hat D)\mathrm d \varphi\mathrm d\xi.\label{equ-M4}
        \end{align}

        \section{The main result}\label{adaptive}
        To estimate the unknown parameter $D$, we construct the following update law
        \begin{align}
                \label{equ-law1}
                \dot{\hat D}=\varrho\mathrm{Proj}_{[\underline D, \overline D]}\{\tau(t)\}, ~~~~0<\varrho<1, 
        \end{align}
        where $\tau(t)$ is given as 
        \begin{align}
                \label{equ-tau}
                \tau(t)=-2\int_0^1\int_{-\pi}^{\pi}(1+s)h(s, \theta, t) P_1(s, \theta, t)\mathrm d\theta\mathrm ds, 
        \end{align}
        and  the standard projection operator is defined as follows
        \begin{align}
                \mathrm{Proj}_{[\underline D,\overline D]}\{\tau(t)\}=\left\{
                \begin{array}{rcl}
                        0 ~~~~    &    & {\hat D=\underline D~and~\tau(t)<0},\\
                        0   ~~~~  &    & {\hat D=\overline D ~and~\tau(t)>0},\\
                        \tau(t)~~     &    & {otherwise}.
                \end{array} \right.\label{equ-law2}
        \end{align} 
        Our claim is that the time-delayed multi-agent system studied in this paper achieves stable formation control, in other words, the state of the error system  \eqref{equ-phi}--\eqref{initial-vartheta} tends to zero under the effect of the adaptive controller \eqref{equ-endU}. The following theorem is established.
        
        \newtheorem{theorem}{\textbf{Theorem}}
        \begin{theorem} \label{theorem1}
                \rm{Consider the closed-loop system consisting of the plant  \eqref{equ-phi}--\eqref{initial-vartheta}, the control law \eqref{equ-endU}, the updated law \eqref{equ-law1}, \eqref{equ-tau} under  Assumption \ref{as1}.  Local boundedness and  \textcolor{black}{regulation} of the system trajectories are guaranteed, i.e., there  exist  positive constants $\mathcal M_1$, $\mathcal{R}_1$ such that if the initial conditions $(\phi_0,\vartheta_0,\hat D_0)$  satisfy $\Psi_{1}(0)<\mathcal M_1$,
                        where 
                        \begin{align}\label{psi}
                                \nonumber \Psi_{1}&(t)=\rVert \phi\rVert^2_{H^2}+\rVert\partial_t \phi\rVert^2_{H^1}+\rVert \vartheta\rVert^2_{H^2}+\rVert \partial_{s\theta\theta } \vartheta\rVert^2+\rVert\partial_{ss\theta} \vartheta\rVert^2\\
                                \nonumber&+\rVert\vartheta(0, \cdot, t)\rVert^2+\rVert\partial_\theta \vartheta(0, \cdot, t)\rVert^2+\rVert\partial_\theta^2 \vartheta(0, \cdot, t)\rVert^2\\
                                &+\rVert\partial_t \vartheta(0, \cdot, t)\rVert^2+\rVert\partial_{t\theta} \vartheta (0, \cdot, t)\rVert^{2}+\tilde D^2,
                        \end{align}
                        the following holds:
                        \begin{align}
                                \Psi_{1}(t)\leq {\mathcal R}_1\Psi_{1}(0), \quad \forall t\geq0\label{equ-Psi};
                        \end{align}
                        furthermore, 
                        \begin{align}
                                &\lim_{t\to \infty}\max_{(s,\theta)\in[0, 1]\times[-\pi,\pi]}|\phi(s, \theta, t)|=0, \label{equ-the-6-phi}\\
                                &\lim_{t\to \infty}\max_{(s,\theta)\in[0, 1]\times[-\pi,\pi]}|\vartheta(s, \theta, t)|=0.\label{equ-the-6-v}
                \end{align}}
        \end{theorem}
        
        \begin{remm} \rm{Only local stability result is obtained due to the existence of the unbounded boundary input operator  combined with the presence of highly nonlinear terms in the target system \eqref{equ-aw0-adp}--\eqref{initial-ah-adp}.
                        In comparison to \cite{WANG2021109909},  the need to ensure continuity of the communication topology of the multi-agent system in three-dimensional space leads to consider more complex norms of the system state (see. \eqref{psi}) for the stability analysis.}\end{remm}

        \section{Proof of the main result}\label{5-proof}
        
        We introduce the following change of variables
        \begin{align}
                m(s, \theta, t)=w(s, \theta, t)-sh(0, \theta, t), \label{equ-m-w}
        \end{align}
        to create a homogeneous boundary condition of the target system
        \begin{align}
                \label{equ-am0-adp}
                &\partial_t{m}(s, {\theta}, t)={\Delta}{m}(s, {\theta}, t)+s\partial_{\theta}^2h(0, \theta, t)-s\partial_t h(0, \theta, t), \\
                &{m}(s, -\pi, t)={m}(s, {\pi}, t),  \quad{m}(0, {\theta}, t)={m}(1, {\theta}, t)=0, \label{bnd-am-adp}\\
                &D\partial_t{h}(s, {\theta}, t)=\partial_s{h}(s, {\theta}, t)-\tilde DP_1(s, \theta, t)-D\dot{\hat D}P_2(s, \theta, t), 
                \label{equ-amh-adp}\\
                &\textcolor{black}{{h}(s, -\pi, t)={h}(s, {\pi}, t)},\quad h(1, {\theta}, t)=0,\label{initial-amh-adp} 
        \end{align}
        with $w(s, \theta, t)$ in $P_{i}(s, \theta, t),~\{ i=1, 2\}$, is rewritten as $m(s, \theta, t)+sh(0, \theta, t)$.
        
        We  will prove Theorem \ref{theorem1} by
        \begin{enumerate}
                \item proving the norm equivalence between the target system \eqref{equ-am0-adp}--\eqref{initial-amh-adp} and the error system  \eqref{equ-phi}--\eqref{initial-vartheta}  through Proposition \ref{proposition4-1},
                \item analyzing  the local stability of the target system \eqref{equ-am0-adp}--\eqref{initial-amh-adp}, and  then deriving the stability of the error system  based on norm equivalence's argument,
                \item and establishing the regulation of the state $\phi(s, \theta, t)$ and $\vartheta(s, \theta, t)$.
        \end{enumerate}
        
        \noindent{\bf (1) Norm equivalence}
        
        We prove the equivalence between the error system  \eqref{equ-phi}--\eqref{initial-vartheta}  and target system  \eqref{equ-am0-adp}--\eqref{initial-amh-adp} in the following Proposition.
        
        \newtheorem{proposition}{\textbf{Proposition}}
        \begin{proposition}\label{proposition4-1}
                \rm{The following estimates hold  between the state of the error system  \eqref{equ-phi}--\eqref{initial-vartheta}, and the state of the target system \eqref{equ-am0-adp}--\eqref{initial-amh-adp}:
                        \begin{align}\label{theo1}
                                \nonumber&\textcolor{black}{\rVert \phi\rVert^2_{H^2}+\rVert\partial_t \phi\rVert^2_{H^1}+\rVert \vartheta\rVert^2_{H^2}}+\rVert \partial_{s\theta\theta } \vartheta\rVert^2+\rVert\partial_{ss\theta} \vartheta\rVert^2\\
                                \nonumber&+\rVert\vartheta(0, \cdot,t)\rVert^2+\rVert\partial_\theta \vartheta(0, \cdot, t)\rVert^2+\rVert\partial_\theta^2 \vartheta(0,\cdot, t)\rVert^2\\
                                &\nonumber+\rVert\partial_t \vartheta(0,\cdot, t)\rVert^2+\rVert\partial_{t\theta} \vartheta (0, \cdot, t)\rVert^2\\
                                \nonumber\leq& R_1(\textcolor{black}{\rVert m\rVert^2_{H^2}+\rVert\partial_tm\rVert^2_{H^1}+\rVert h\rVert^2_{H^2}}+\rVert \partial_{s\theta\theta }h\rVert^2+\rVert\partial_{ss\theta}h\rVert^2\\
                                \nonumber&+\rVert h(0, \cdot, t)\rVert^2+\rVert\partial_\theta h(0, \cdot, t)\rVert^2+\rVert\partial_\theta^2h(0,  \cdot, t)\rVert^2\\
                                &+\rVert\partial_th(0,  \cdot, t)\rVert^2+\rVert\partial_{t\theta}h (0,  \cdot, t)\rVert^{2}),\\
                                \nonumber&\textcolor{black}{\rVert m\rVert^2_{H^2}+\rVert\partial_tm\rVert^2_{H^1}+\rVert h\rVert^2_{H^2}}+\rVert \partial_{s\theta\theta }h\rVert^2+\rVert\partial_{ss\theta}h\rVert^2\\
                                \nonumber&+\rVert h(0,  \cdot, t)\rVert^2+\rVert\partial_\theta h(0,  \cdot, t)\rVert^2+\rVert\partial_\theta^2h(0,  \cdot, t)\rVert^2\\
                                \nonumber&+\rVert\partial_th(0,  \cdot, t)\rVert^2+\rVert\partial_{t\theta}h (0, \cdot, t)\rVert^{2})\\
                                \nonumber\leq&R_2(\textcolor{black}{\rVert \phi\rVert^2_{H^2}+\rVert\partial_t \phi\rVert^2_{H^1}+\rVert \vartheta\rVert^2_{H^2}}+\rVert \partial_{s\theta\theta } \vartheta\rVert^2+\rVert\partial_{ss\theta} \vartheta\rVert^2\\
                                \nonumber&+\rVert\vartheta(0, \cdot, t)\rVert^2+\rVert\partial_\theta \vartheta(0,  \cdot, t)\rVert^2+\rVert\partial_\theta^2 \vartheta(0,  \cdot, t)\rVert^2\\
                                &+\rVert\partial_t \vartheta(0,  \cdot, t)\rVert^2+\rVert\partial_{t\theta} \vartheta (0,  \cdot, t)\rVert^2),\label{theo2}
                        \end{align}
                        where $R_i$,  $i=1, 2$ are sufficiently large positive constants.}
        \end{proposition} 
        
       \textcolor{black}{The proof of Proposition \ref{proposition4-1} is stated in Appendix \ref{apA}. }
        
        Next, we  show the local stability for the closed-loop system consisting of the  $(\phi,  \vartheta)$-system under the control law \eqref{equ-endU}, and  with the updated law \eqref{equ-law1}--\eqref{equ-tau}.

        \noindent{\bf (2) Local stability analysis}

        Since the error system \eqref{equ-phi}--\eqref{initial-vartheta} is equivalent to the target system \eqref{equ-am0-adp}--\eqref{initial-amh-adp},  we  establish the local stability of  the target system by introducing the following  Lyapunov-Krasovskii-type function,      
\begin{align}
\label{equ-6-V-adp}
\nonumber V_{1}&(t)=b_1 (\textcolor{black}{\|m\|^2_{H^2}+\|\partial_tm\|^2_{H^1}})+D\int_0^1\int_{-\pi}^\pi(1+s)(|h|^2\\
&\nonumber+|\partial_sh|^2+|\partial_\theta h|^2+|\Delta h|^2+|\partial_{s\theta\theta}h|^2+|\partial_{ss\theta}h|^2)\mathrm d\theta\mathrm ds\\ 
&\nonumber+b_2D(\textcolor{black}{\|h(0, \cdot, t)\|^2+\|\partial_\theta h(0, \cdot, t)\|^2+\|\partial_\theta^2 h(0, \cdot, t)\|^2}\\
&\textcolor{black}{+\|\partial_t h(0, \cdot, t)\|^2+\|\partial_{t\theta}h(0, \cdot, t)\|^2})+\frac{\tilde D^{2}}{2\varrho}.
\end{align}
Taking the time derivative  of \eqref{equ-6-V-adp}, \textcolor{black}{based on \eqref{equ-P1}, \eqref{equ-P2}, \eqref{equ-m-w}--\eqref{initial-amh-adp},}  and using Cauchy Schwartz's inequality, Young's inequality, \textcolor{black}{Poincare's inequality,} and integration by parts, we obtain
we obtain that
        \begin{align}
                \label{equ-P11-adp}
                \nonumber\dot{V}_{1}&(t)\leq-b_1\bigg(\frac{3}{8}-\frac{1}{\sigma_2}-\frac{1}{\sigma_3}\bigg)\rVert{m}\rVert^2-2b_1\rVert\partial_\theta{m}\rVert^2\\
                \nonumber&-\frac{b_1}{2}\rVert\partial_s{m}\rVert^2-b_1\bigg(2-\frac{1}{\sigma_{1}}-\frac{1}{\sigma_{4}}-\frac{1}{\sigma_{5}}\bigg)\rVert\Delta m\rVert^2-b_1\bigg(\frac{3}{8}\\
                \nonumber&-\frac{1}{\sigma_6}-\frac{1}{\sigma_7}\bigg)\rVert\partial_t{m}\rVert^2-\frac{b_1}{2}\rVert\partial_{ts}{m}\rVert^2-2b_1\rVert\partial_{t\theta}{m}\Vert^2-b_1(2\\
                \nonumber&-\sigma_{1}-\sigma_{8}-\sigma_{9})\rVert\Delta\partial_t {m}\rVert^2-(\rVert h\rVert^2+\rVert\partial_sh\rVert^2+\rVert\partial_\theta h\rVert^2\\
                \nonumber&+\rVert\Delta h\rVert^2+\rVert\partial_{s\theta  \theta }h\rVert^2+\rVert\partial_{ss\theta}h\rVert^2)-(1-b_{2}\sigma_{10})\textcolor{black}{\rVert h(0, \cdot, t)\rVert^2}\\
                \nonumber&-\bigg(1-\frac{b_{2}}{\sigma_{10}}-\frac{3b_2\sigma_{13}}{D}-\frac{b_1(\sigma_3+\sigma_5)}{D^{2}}\bigg)\textcolor{black}{\rVert \partial_sh(0, \cdot, t)\rVert^2}\\
                \nonumber&-\bigg(1-\frac{b_{2}}{\sigma_{11}}\bigg)\textcolor{black}{\rVert \partial_\theta h(0, \cdot, t)\rVert^2}-\bigg(1-\frac{7b_2}{D^3\sigma_{13}}-\frac{7b_1}{3D^4}({\sigma_7}\\
                \nonumber&+\frac{1}{\sigma_9})\bigg)\textcolor{black}{\rVert\partial_s^2h(0, \cdot, t)\rVert^2}-\bigg(1-\frac{b_1(\sigma_2+\sigma_{4})}{3}-\frac{b_2}{\sigma_{12}}\bigg)\\
                \nonumber&\cdot\textcolor{black}{\rVert \partial^2_{\theta}h(0, \cdot, t)\rVert^2}-\bigg(1-\frac{7b_{2}}{D^3\sigma_{14}}\bigg)\textcolor{black}{\rVert\partial_{ss\theta}h(0, \cdot, t)\rVert^2}-\bigg(2\\
\nonumber&-b_{2}\sigma_{11}-\frac{3b_{2}\sigma_{14}}{D}\bigg)\textcolor{black}{\rVert\partial_{s\theta}h(0, \cdot, t)\rVert^2}-\bigg(1-\sigma_{12}b_{2}\\
                \nonumber&-\frac{b_1}{D^2}({\sigma_6}+\frac{1}{\sigma_{8}})\bigg)\textcolor{black}{\rVert\partial_{s\theta\theta}h(0, \cdot, t)\rVert^2}-\tilde DE_{1}(t)-D\dot{\hat D}E_{2}(t)\\
&+\tilde D^2E_{3}(t)+\dot{\hat D}^2E_{4}(t)+\ddot{\hat D}^2E_{5}(t)-\dot{\hat D}\frac{\tilde D}{\varrho}, 
        \end{align}
        where $\sigma_i>0$, $i=1,2,...,14$, and 
        \begin{align}
                \nonumber E_1&(t)=2\int_{0}^{1}\int_{-\pi}^{\pi}(1+s)(hP_1+\partial_sh\partial_sP_1+\partial_\theta h\partial_\theta P_1\\
                \nonumber&+\Delta h\partial^2_s P_1+\Delta h\partial^2_\theta P_1+\partial _{s\theta\theta} h\partial_{s\theta\theta}P_1+\partial_{ss\theta} h\partial_{ss\theta} P_1)\mathrm{d}\theta\mathrm{d}s\\
                \nonumber&+2b_2\int_{-\pi}^{\pi}(h(0, \theta, t)P_1(0, \theta, t)+\partial_\theta h(0, \theta, t)\partial_\theta P_1(0, \theta, t)\\
                &+\partial_\theta^2 h(0, \theta, t)\partial_\theta^2 P_1(0, \theta, t)){\mathrm{d}\theta}{\mathrm{d}s}, \\
                \nonumber E_2&(t)=2\int_{0}^{1}\int_{-\pi}^{\pi}(1+s)(hP_2+\partial_sh\partial_sP_2+\partial_\theta h\partial_\theta P_2\\
                \nonumber&+\Delta h\partial^2_s P_2+\Delta h\partial^2_\theta P_2+\partial_{s\theta\theta} h\partial_{s\theta \theta} P_2+\partial_{ss\theta } h\partial_{ss\theta } P_2){\mathrm{d}\theta}{\mathrm{d}s}\\
                \nonumber&+2b_2\int_{-\pi}^{\pi}(h(0, \theta, t)P_2(0, \theta, t)+\partial_\theta h(0, \theta, t)\partial_\theta P_2(0, \theta, t)\\
                &+\partial_\theta^2 h(0, \theta, t)\partial_\theta^2 P_2(0, \theta, t)){\mathrm{d}\theta}{\mathrm{d}s}, \\
                \nonumber E_3&(t)=\bigg(\frac{b_1(\sigma_3+\sigma_5)}{D^2}+\frac{3b_{2}\sigma_{13}}{D}\bigg)\rVert  P_1(0, \cdot, t)\rVert^{2}+\frac{b_1}{D^2}\bigg({\sigma_6}\\
                \nonumber&+\frac{1}{\sigma_8}\bigg)\rVert \partial_\theta^2 P_1(0, \cdot, t)\rVert ^{2}+\frac{3b_{2}\sigma_{14}}{D}\rVert \partial_\theta P_1(0, \cdot, t)\rVert ^{2}\\
                \nonumber&+{7}\bigg(\frac{b_1}{3D^{2}}({\sigma_{7}}+\frac{1}{\sigma_{9}})+\frac{b_{2}}{\sigma_{14}D}\bigg)\bigg(\frac{1}{D^2}\rVert \partial_sP_1(0, \cdot, t)\rVert ^{2}\\
                \nonumber&+\rVert \partial_t P_1(0, \cdot, t)\rVert ^{2}\bigg)+\frac{7b_{2}}{\sigma_{15}D}\bigg(\frac{1}{D^2}\rVert \partial_{s\theta}P_1(0, \cdot, t)\rVert ^{2}\\
                \nonumber&+\rVert \partial_{t\theta} P_1(0, \cdot, t)\rVert ^{2}\bigg)+4(\rVert P_1(1,\cdot,t)\rVert ^2+2\rVert \partial_{\theta} P_1(1, \cdot, t)\rVert ^{2}\\
                \nonumber&+\rVert \partial_{\theta}^2 P_1(1, \cdot, t)\rVert ^{2})+12(\rVert \partial_{s}P_1(1, \cdot, t)\rVert ^{2}+\rVert \partial_{s\theta}P_1(1, \cdot, t)\rVert ^{2}\\
                &+D^{2}\rVert \partial_t P_1(1, \cdot, t)\rVert ^{2}+D^{2}\rVert \partial_{t\theta}P_1(1, \cdot, t)\rVert ^{2}), \\
                \nonumber E_4&(t)=\bigg(\frac{7b_1}{3D^{2}}({\sigma_{7}}+\frac{1}{\sigma_{9}})+\frac{7b_{2}}{D\sigma_{13}}\bigg)(\rVert P_1(0, \cdot, t)\rVert^{2}\\
                \nonumber&+\rVert\partial_s P_2(0, \cdot, t)\rVert^{2}+D^{2}\rVert\partial_t P_2(0, \cdot, t)\rVert^{2})+b_1\bigg({\sigma_6}+\frac{1}{\sigma_8}\bigg)\\
                \nonumber&\cdot\rVert\partial_\theta^2 P_2(0, \cdot, t)\rVert^{2}+\frac{7b_{2}}{D\sigma_{14}}(\rVert\partial_\theta P_1(0, \cdot, t)\rVert^{2}+\rVert\partial_{s\theta} P_2(0, \cdot, t)\rVert^{2}\\
                \nonumber&+D^{2}\rVert\partial_{t\theta} P_2(0, \cdot, t)\rVert^{2})+3b_{2}D\sigma_{14}\rVert\partial_\theta P_2(0,\cdot, t)\rVert^{2}+(b_1(\sigma_{3}\\
                \nonumber&+\sigma_{5})+3b_{2}D\sigma_{13})\rVert P_2(0,\cdot, t)\rVert^{2}+4D^{2}(\rVert P_2(1, \cdot, t)\rVert^{2}\\
                \nonumber&+2\rVert\partial_{\theta} P_2(1, \cdot, t)\rVert^{2}+\rVert\partial_{\theta}^2 P_2(1, \cdot, t)\rVert^{2})+12D^{2}(\rVert P_1(1, \cdot, t)\rVert^{2}\\
                \nonumber&+\rVert\partial_{s}P_2(1, \cdot, t)\rVert^{2}+\rVert\partial_\theta P_{1}(1, \cdot, t)\rVert^{2}+D^{2}\rVert\partial_tP_2(1, \cdot, t)\rVert^{2}\\
&+\rVert\partial_{s\theta}P_2(1,\cdot, t)\rVert^{2}+D^{2}\rVert\partial_{t\theta} P_2(1, \cdot, t)\rVert^{2}), \\
                \nonumber E_5&(t)=\bigg(\frac{7b_1}{3}(\sigma_{7}+\frac{1}{\sigma_{9}})+\frac{7b_{2}D}{\sigma_{13}}\bigg)\rVert P_2(0, \cdot, t)\rVert^{2}+\frac{7Db_{2}}{\sigma_{14}}\\
                &\cdot\rVert \partial_{\theta} P_2(0, \cdot, t)\rVert^{2}+12D^{4}(\rVert P_2(1, \cdot, t)\rVert^{2}+\rVert\partial_{\theta}P_2(1, \cdot, t)\rVert^{2}).
        \end{align}
        By setting $\sigma_1=1$, $\sigma_2=\sigma_3=8$, $\sigma_4=\sigma_5=3$, $\sigma_6=\sigma_7=8$, $\sigma_8=\sigma_9=\frac{1}{3}$, $\sigma_{10}=\sigma_{11}=\sigma_{12}=\sigma_{13}=\sigma_{14}=1$, $0<b_1<\min\{\frac{3}{11}$, $\frac{\underline D^2}{11}$, $\frac{3\underline D^4}{77}\}$, $0<b_2<\min\{\frac{2\underline D}{3+\overline D}$, $ \frac{3-11b_{1}}{3}$, $\frac{\underline D^2-11b_{1}}{\overline D^2}$, $\frac{\underline D^3}{7}$, $\frac{\underline D^2-11b_1}{\overline D(3+2\overline D)}$, $\frac{3\underline D^4-77b_1}{3\overline D(\overline D^3+7)}\}$, we get the following estimate 
        \begin{align}
                \label{equ-P12-adp}
                \nonumber\dot{V}_1(t)\leq&-\kappa_{1} V_2(t)-\tilde DE_{1}(t)-D\dot{\hat D}E_{2}(t)+\tilde D^2E_{3}(t)\\
                &+\dot{\hat D}^2E_{4}(t)+\ddot{\hat D}^2E_{5}(t)-\dot{\hat D}\frac{\tilde D}{\varrho}, 
        \end{align}
        where $\kappa_1=\min\{\frac{b_1}{8},~1-\frac{11b_1}{\underline D^{2}}-{b_{2}}-\frac{3b_{2}}{\underline D}\}> 0$ and 
        \begin{align}\label{equ-V0-ada}
                \nonumber V_2(t)=&\rVert m\rVert^2_{H^2}+\rVert\partial_tm\rVert^2_{H^1}+\rVert h\rVert^2_{H^2}+\rVert \partial_{s\theta\theta}h\rVert ^2+\rVert \partial_{ss\theta}h\rVert ^2\\
                &+\rVert h(0, \cdot, t)\rVert^2.
        \end{align}
        With the help of Agmon's, Cauchy-Schwarz, and Young's inequalities, one can perform quite long calculations to derive 
        the following estimates:
        \begin{align}\label{equ-V0-L0}
                &E_1(t)\leq11L_{1}V_2(t), \quad\quad E_2(t)\leq11L_{1}V_2(t), \\
                &E_3(t)\leq(\alpha_1+\alpha_2\theta^2 L_{1}^2V_2(t)^2+\tilde D^2\alpha_2)L_{1}V_2(t), \\
                &E_4(t)\leq (\alpha_3+\alpha_4\theta^2L_{1}^2V_4(t)^2+\tilde D^2\alpha_4)L_{1}V_2(t), \\
                &E_5(t)\leq \alpha_4L_{1}V_2(t), \quad \dot{\hat D}(t)\leq \theta L_{1}V_2(t), \\
                &\ddot{\hat D}(t)\leq\theta(1+\theta L_{1}V_2(t)+|\tilde{D}|)L_{1}V_2(t), \label{equ-V0-L1}
        \end{align}
        where
        \begin{align}
                \alpha_1&=\frac{11(13\overline D^{2}+7)b_1}{3\underline D^{4}}+\frac{2(10\overline D^{2}+7)b_{2}}{\underline D^3}+24\overline D^{2}+40,\\
                \alpha_2&=\frac{77b_1}{3\underline D^2}+\frac{14b_{2}}{\underline D}+24\overline D^{2}, \\
                \alpha_3&=\frac{11(13\overline D^{2}+14)b_1}{3\underline D^{2}}+\frac{4(5\overline D^{2}+7)b_{2}}{\underline D}+64\overline D^{2}+24\overline D^4\\
                \alpha_4&=\frac{77b_1}{3}+{14\overline Db_{2}}+24\overline D^{4}, 
        \end{align}
        and $L_1$ is a sufficiently large positive constant, which estimation method is similar to the method in Appendix \ref{apA}. And then, combining with \eqref{equ-V0-L0}--\eqref{equ-V0-L1}, one can get
\begin{align}
\label{equ-6-dot_V}
\nonumber\dot{V}_{1}&(t)\leq-\kappa_1 V_2(t)+|\tilde D|(8+\frac{1}{\varrho})L_{1}V_{2}(t)+8\overline DL_{1}^2V_{2}(t)^2\\
\nonumber&+\tilde D^2\alpha_1L_{1}V_2(t)+\tilde D^2(2\alpha_2+12\alpha_4)L_{1}^3V_2(t)^3+\tilde D^4\alpha_2L_{1}\\
&\cdot V_2(t)+(\alpha_3+12\alpha_4)L_{1}^3V_2(t)^3+28\alpha_4L_{1}^5V_2(t)^5.
\end{align}
From \eqref{equ-6-V-adp}, it is easy to get $\tilde{D}^2\leq 2\varrho V_{1}(t)-2\varrho\zeta_1V_2(t)$,
$\zeta_1=\min\{b_1, ~\underline D, ~b_2\underline D\}$. Using Cauchy-Schwarz's and Young's inequalities, one can deduce that
\begin{align}\label{equ-tilde-D}
\left|\tilde{D}\right|\leq& \frac{\varepsilon_1}{2}+\frac{\tilde D^2}{2\varepsilon_1}\leq\frac{\varepsilon_1}{2}+\frac{\varrho}{\varepsilon_1}V_1(t)-\frac{\varrho\zeta_1}{\varepsilon_1} V_2(t). 
\end{align}
%
Again, using \eqref{equ-6-V-adp} we have
\begin{align}
\zeta_1 V_0(t)\leq V(t).\label{VI}
\end{align}
\textcolor{black}{Substituting \eqref{equ-tilde-D}, \eqref{VI} into \eqref{equ-6-dot_V}, we derive the following estimate}
\begin{align}
\nonumber\dot V_{1}&(t)\leq-\bigg(\frac{\kappa_1}{2} -8\varrho^2\alpha_2 L_{1}V_1(t)^2\bigg)V_2(t)-\bigg(\frac{\kappa_1}{2} -L_{1}(8\\
\nonumber&+\frac{1}{\varrho})(\frac{\varepsilon_1}{2}+\frac{\varrho}{\varepsilon_1}V_{1}(t))-2\varrho\alpha_1 L_1V_1(t)\bigg)V_2(t)-L_{1}\\
\nonumber&\cdot\bigg(\frac{\varrho\zeta_1 }{\varepsilon_1} (8+\frac{1}{\varrho})-(\frac{(\alpha_3+12\alpha_4)L_{1}^{2}}{\zeta_1 } +8\alpha_2\varrho^2\zeta_1 ^{2})V_1(t)\\
\nonumber&-8\overline DL_{1}\bigg)V_2(t)^2-2{\varrho }L_1\bigg({\alpha_1}\zeta_1-\frac{(\alpha_2+13\alpha_4)L_{1}^2}{\zeta_1 } V_1(t)^2\bigg)  \\
&\cdot V_2(t)^2-L_{1}^3\bigg(2\varrho(\alpha_2+13\alpha_4)\zeta_1-\frac{28\alpha_4L_{1}^{2}}{\zeta_1 } V_1(t)\bigg) V_2(t)^4.\label{equ-V1-last}
\end{align}
Let $\varepsilon_1$ defined  as $ \varepsilon_1<\min\left\{\frac{\kappa_1\varrho}{L_{1}(8\varrho+1)}, \frac{(8\varrho+1)\zeta_1 }{8\varrho\overline D L_{1}}\right\}$,
to ensure $V_{1}(0)\leq \mathcal \mu_1$, where
\begin{align}
\label{equ-initial-coditions1}
 \nonumber\mu_1\triangleq&\min\left\{\frac{\varepsilon_1(\kappa_1{\varrho }-(8\varrho+1)L_{1}\varepsilon_1)}{2\varrho L_{1}(8\varrho+1+2\alpha_1\varrho\varepsilon_1))}, {\frac{\sqrt{\kappa_1} }{4\varrho\sqrt{\alpha_2L_{1}}} },\right.\\
\nonumber& \frac{\sqrt{\alpha_1}\zeta_1  }{\sqrt{(\alpha_2+13\alpha_4)}L_{1}},\frac{\varrho(\alpha_2+13\alpha_4)\zeta_1  }{14\alpha_4L_{1}^{2}},\\
&\left.\frac{\zeta_1 ((8\varrho+1)\zeta_1 -8\overline DL_{1}\varepsilon_1) }{\varepsilon_1(4\varrho^{2}\alpha_2\zeta_1^2+(\alpha_3+12\alpha_4)L_{1}^{2})}\right\}.
\end{align}
Therefore,
\begin{align}
\nonumber\dot V_{1}(t)\leq&-(\delta_1(t)+\delta_2(t))V_2(t)-(\delta_3(t)+\delta_4(t))V_2(t)^2\\
&+\delta_5(t)V_2(t)^4,
\end{align}
where
\begin{align}
\delta_1(t)=&\frac{\kappa_1}{2} -L_{1}(8+\frac{1}{\varrho})(\frac{\varepsilon_1}{2}+\frac{\varrho}{\varepsilon_1}V_{1}(t))-2\varrho\alpha_1 L_{1}V_1(t), \label{equ-6-delta1}\\
\delta_2(t)=&\frac{\kappa_1}{2} -4\varrho^2\alpha_2 L_{1}V_1(t)^2, \\
\nonumber\delta_3(t)=&L_{1}\bigg(\frac{\varrho\zeta_1  }{\varepsilon_1}(8+\frac{1}{\varrho}) -8\overline DL_{1}-(\frac{(\alpha_3+12\alpha_4)L_{1}^{2}}{\zeta_1 } \\
&+4\varrho^2\alpha_2)V_1(t)\bigg), \\
\delta_4(t)=&2{\varrho}L_{1}\left({\alpha_1}\zeta_1 -\frac{(\alpha_2+13\alpha_4)L_{1}^{2}}{\zeta_1 } V_1(t)^2\right), \\
\delta_5(t)=&L_{1}^3\bigg({2{\varrho }{(\alpha_2+13\alpha_4)}}\zeta_1 -\frac{28\alpha_4L_{1}^{2}}{\zeta_1 } V_1(t)\bigg),
\end{align}
are nonnegative functions if the initial condition satisfies \eqref{equ-initial-coditions1}. Thus, $V_{1}(t)\leq V_{1}(0),~\forall t\geq0$.
        
Using \eqref{theo1}, we can get
\begin{align}
 \Psi_{1}(t)
\leq& \frac{\max\{R_1, 1\}}{\min\{b_1,  2\overline D, b_2\overline D, \frac{1}{2\varrho}\}}V_{1}(t)\leq{\mu}_2V_{1}(0), \label{equ-initial-coditions2}
\end{align}
where $\Psi_{1}(t)$ is defined as \eqref{psi}, and $\mu_2=\frac{\max\{R_1, 1\}}{\min\{b_1,~ \underline D, ~b_2\underline D, ~\frac{1}{2\varrho}\}}$.
Hence, combining \eqref{equ-initial-coditions1} and \eqref{equ-initial-coditions2}, we have $\mathcal M_1=\mu_1\mu_2$.
        
From \eqref{theo2} and \eqref{equ-6-V-adp}, one gets
        \begin{align}
                \nonumber V_{1}(t)
                \leq&\max\left\{\max\{b_1,  2\overline D, b_2\overline D\}R_{2}, \frac{1}{2\varrho}\right\}\Psi_{1}(t).
        \end{align}
        Knowing that $V_{1}(0)\leq\max\{\max\{b_1,  2\overline D, b_2\overline D\}R_{2}, \frac{1}{2\varrho}\}\Psi_{1}(0)$,
        we arrive at \eqref{equ-Psi}  with ${\mathcal{R}}_{1}=\mu_2\max\{\max\{b_1,  2\overline D, b_2\overline D\}R_{2}, \\\frac{1}{2\varrho}\}$, 
        which proves the local stability of the closed-loop system.
        
        Next, we will prove the regulation of the cascaded system $(\phi,\vartheta)$ to complete the  proof of Theorem 1.
        
        \noindent{\bf (3) Regulation of the cascaded system}
        
        From \eqref{equ-6-V-adp} and \eqref{equ-V1-last}, we get the boundedness of all terms in  \eqref{equ-V0-ada}, and then, based on \eqref{theo1}, we also get the boundedness of all terms of  $\Psi_{1}(t)$. We will prove \eqref{equ-the-6-phi} and \eqref{equ-the-6-v} in Theorem \ref{theorem1} by applying  Lemma D.2 \cite{Krstic2010} to ensure the following facts:
        
        \begin{itemize}
                \item  all terms  in \eqref{equ-V0-ada} are square integrable in time,
                \item $\frac{\mathrm d}{\mathrm dt}(\rVert m\rVert^2)$,  $\frac{\mathrm d}{\mathrm dt}(\rVert h\rVert^2)$ and $\frac{\mathrm d}{\mathrm dt}(\rVert\partial_{s}h\rVert^2)$ are bounded.
        \end{itemize}
        
        Knowing that
        \begin{align}
                \int_0^t\rVert m(\tau)\rVert^2\mathrm d\tau\leq\frac{1}{\inf_{0\leq\tau\leq t}\delta_1(\tau)}\int_0^t\delta_1(\tau)V_2(\tau)\mathrm d\tau, \label{equ-square-w}
        \end{align} 
        and using  \eqref{equ-6-delta1}, the following inequality holds:
        \begin{align}
                \nonumber \inf_{0\leq\tau \leq t}\delta_1(t)=&\frac{\kappa_1}{2} -L_{1}(8+\frac{1}{\varrho})\left(\frac{\varepsilon_1}{2}
                +\frac{\varrho}{\varepsilon_1}V_{1}(t)\right)\\
                &-2\varrho\alpha_1 L_{1}V_1(t).\label{equ-square1}
        \end{align} 
        Since $\dot V_{1}\leq-(\delta_1(t)+\delta_2(t))V_2(t)-(\delta_3(t)+\delta_4(t))V_1(t)^2+\delta_5(t)
        V_2(t)^4$ and $\delta_i(t)$ are nonnegative  functions, we have $\dot V_1\leq-\delta_1(t)V_2(t)$, and 
       integrating it over $[0,\ t]$ leads to 
        \begin{align}
                \int_0^t\delta_1(\tau)V_2(\tau)\mathrm d\tau\leq V_1(0)\leq \mu_1.\label{equ-square2}
        \end{align}
        Substituting \eqref{equ-square1} and \eqref{equ-square2} into \eqref{equ-square-w}, we get $\rVert m\rVert$ is square integrable in time. Similarly, one can establish that other terms in   \eqref{equ-V0-ada} are square-integrable in time. 
        
        To prove that  $\frac{\mathrm d}{\mathrm dt}(\rVert m\rVert^2)$,  $\frac{\mathrm d}{\mathrm dt}(\rVert h\rVert^2)$ and $\frac{\mathrm d}{\mathrm dt}(\rVert\partial_{s}h\rVert^2)$ are bounded, we define the Lyapunov function
        \begin{align}
 &V_3(t)=\frac{1}{2}\textcolor{black}{\| m\|^2}+\frac{b_3D}{2}\int_0^1\int_{-\pi}^\pi(1+s)(|h|^2+|\partial_sh|^2)\mathrm d\theta\mathrm ds, \label{equ-Vnew3}
        \end{align}
        where $b_3$ is a positive constant. Taking the derivative of \eqref{equ-Vnew3} with respect to time, 
        and using integration by parts and Young's inequality, the following holds
        \begin{align}
                \nonumber\dot V_3&(t)\leq-\rVert \partial_sm\rVert^2-b_3\rVert h\rVert^2-b_3\rVert \partial_sh\rVert^2+(\frac{1}{2\iota_7}+\frac{1}{2\iota_8})\rVert m\rVert^2\\
                \nonumber&+\frac{\iota_7}{6}\rVert \partial_{\theta\theta}h\rVert^2+\frac{\iota_7}{6}\rVert \partial_{\theta\theta s}h\rVert^2+\frac{\iota_8}{2\underline D^2}|\tilde D|^2\rVert P_1(0, \cdot, t)\rVert^2\\
                \nonumber&-(\frac{b_3}{2}-\frac{\iota_8}{2\underline D^2})\rVert\partial_sh(0, \cdot, t)\rVert^2+\frac{\iota_8}{2}|\dot{\hat D}|^2\rVert P_2(0, \cdot, t)\rVert^2+4b_3\\
                \nonumber&\cdot|\tilde D|^{2}\rVert h\rVert\rVert P_1\rVert+2b_3|\tilde D|^2\rVert P_1(1, \cdot, t)\rVert^2+4b_3{|\tilde D|}\rVert h\rVert \rVert P_{1}\rVert\\
                \nonumber&+2b_3\overline D^{2}|\dot{\hat D}|^2\rVert P_2(1, \cdot, t)\rVert^2+4b_3|\tilde D|^{2}\rVert  \partial_sh\rVert\rVert \partial_sP_1\rVert\\
                &+4b_3{|\tilde D|}\rVert h\rVert \rVert P_{2}\rVert+4b_3\overline D|\dot{\hat{D}}|\rVert \partial_sh\rVert \rVert\partial_sP_{2}\rVert.
        \end{align}
        Setting $\iota_7=\iota_8=8$ and $b_3>\frac{4}{\underline D^2}$, we have \begin{align}
               \dot V_3 \leq&-c_{1}V_3+f_1(t)V_{3}+f_2(t)<\infty, \label{equ-6-wt}
        \end{align}
        where we use Young's and Agmon's inequalities,  $c_{1}=\min\{\frac{1}{4},~\frac{1}{2\overline D}\}$, and
        \begin{align}
                f_1&(t)=\frac{{2}\tilde D^{2}}{\underline D}(\tilde D^{2}+\overline D^2|\dot{\hat{D}}|^2), \\
                \nonumber f_2&(t)=\frac{4}{3}(\rVert \partial_{\theta\theta}h\rVert^2+\rVert \partial_{\theta\theta s}h\rVert^2)+2b_3|\tilde D|^2\rVert P_1(1, \cdot, t)\rVert^2\\
                \nonumber&
                +2b_3\overline D^{2}|\dot{\hat D}|^2\rVert P_2(1, \cdot, t)\rVert^2+\frac{4|\tilde D|^2}{\underline D^2} \rVert P_1(0, \cdot, t)\rVert^2+4|\dot{\hat D}|^2\\
                \nonumber&\cdot\rVert P_2(0, \cdot, t)\rVert^2+2b_3|\tilde D|^2\rVert P_1(1, \cdot, t)\rVert^2+2b_3\rVert P_1\rVert^2\\
                \nonumber&+2b_3\rVert \partial_s P_1\rVert^2+2b_3\overline D^{2}|\dot{\hat D}|^2\rVert P_2(1, \cdot, t)\rVert^2+2b_3\rVert P_2\rVert^2\\
                &+2b_3\rVert \partial_s P_2\rVert^2. 
        \end{align}
        Combining \eqref{equ-P1} and \eqref{equ-P2}, we get that $|\dot{\hat D}|$, $\rVert P_1(0,\cdot, t)\rVert^2$, $\rVert P_2(0,\cdot, t)\rVert^2$, $\rVert P_1(1,\cdot, t)\rVert^2$, $\rVert P_2(1,\cdot, t)\rVert^2$, $\rVert P_1\rVert^2$ and $\rVert P_2\rVert^2$  are  bounded and integrable. Thereby,
        $f_1(t)$ and $f_2(t)$ are bounded and integrable functions of time. Thus, from \eqref{equ-6-wt}, we deduce that $\dot V_2\leq\infty$, which proves the boundedness of $\frac{\mathrm d}{\mathrm dt}(\rVert m\rVert^2)$,  $\frac{\mathrm d}{\mathrm dt}(\rVert h\rVert^2)$ and $\frac{\mathrm d}{\mathrm dt}(\rVert\partial_{s}h\rVert^2)$. Moreover, by Lemma D.2 \cite{Krstic2010}, it holds that   $\rVert m\rVert$, $\rVert h\rVert$, $ \rVert\partial_{s}h\rVert\to 0$ as $t\to\infty$. Knowing that  $\rVert h(0,\cdot,t)\rVert^2\leq 2\rVert h\rVert\rVert\partial_{s}h\rVert$, so $\rVert h(0,\cdot,t)\rVert^2\to 0$ as $t\to\infty$. From \eqref{equ-trans3} and \eqref{equ-m-w}, one can get 
\begin{align}
\rVert \phi\rVert^{2}\leq&4(1+\|l(s, \tau)\|^{2})\rVert m\rVert^{2}+4\|l(s, \tau)\|^{2}(\rVert h\rVert^{2}+\rVert \partial_sh\rVert^2). 
        \end{align}
        So, we get $\rVert\phi\rVert^{2}\to 0$ as $t\to\infty$. Since $\rVert\phi\rVert_{H^2}$ is bounded, we can get $\phi(s, \theta, t)^2\leq C\rVert \phi\rVert\rVert \phi\rVert_{H^2}$  by using Agmon's inequality, and then we get  $\phi(s,\theta,t)$ is regulated. Similarly, we can get $\vartheta(s,\theta,t)$ is also regulated.

        \section{Numerical Simulations}\label{simu}
        \subsection{Control laws for the leaders and the followers} In order to implement  control laws of the followers, we  discretize the  PDEs  \eqref{equ-u_0} and \eqref{equ-z_0}. For $u\in\Omega$, we define the following  discretized grid
        \begin{align}
                \label{equ-disway} 
                s_i=(i-1)h_s, ~~\theta_j=(j-1)h_\theta,~~ d_k=(k-1)\Delta D,
        \end{align}
        for $i=2,...,M-1$, $j=1,...,N$, $k=1,...,M'$, where 
        $h_s=\frac{1}{M-1}$,
        $h_\theta=\frac{2\pi}{N-1}$ and $\Delta D=\frac{D}{M'-1}$.
        Using a three-point central difference approximation, the control laws of the follower agents $(i,j)$ are written as
        \begin{align}
                \label{equ-disu} 
                \dot{u}_{ij}=&\frac{(u_{i+1,j}-u_{i,j})
                        -(u_{i,j}-u_{i-1,j})}{h_s^2}+\beta_1\frac{u_{i+1,j}-u_{i-1,j}}{2h_s}
                \nonumber\\ &+\frac{(u_{i,j+1}-u_{i,j})
                        -(u_{i,j-}u_{i,j-1})}{h_\theta^2}
                +\lambda_1{u}_{i,j},
        \end{align}
        where  $i=2,...,M-1$, $j=1,...,N$,  and  all the state variables in $\theta$ space are $2\pi$ periodic, namely, $u_{i,1}=u_{i,N}$. The leader agents with  guiding role  at the boundary $s=0$, namely  $i=1$ are formed as $u_{1,j}=f_1(\theta_j)$.
        For the leader agents at the boundary $s=0$, namely $i=M$, from the discretized form of \eqref{equ-endU}, the  state feedback control action is given by
        \begin{align}
                \label{equ-disU} 
                u_{M,j}(t)&= \sum_{m=1}^{M}\sum_{l=1}^{N}{a}_{m,l}\gamma_{j,m,l}e^{-\frac{1}{2}\beta_1(1-s_m)}(u_{m,l}(t)-\bar{u}_{m,l}(t))
                \nonumber\\ &-\sum_{k=1}^{M'}\sum_{l=1}^{N}a'_{k,l}
                \gamma'_{j,k,l}u_{M,l}(t-D+d_k)+\overline u_{M,j},
        \end{align}
        where $\gamma_{j,m,l}$ and $\gamma'_{j,k,l}$ can be discretized from \eqref{kernel-gamma} and \eqref{kernel-p}.
        $M$, $N$, and $M'$ are  odd numbers according to Simpson's rule. The control laws for the  $z$-coordinate can be obtained in a similar way.

        
        \subsection{Simulation results}
\begin{figure}[htbp]
                \centering
                \subfigure[]{\includegraphics[width=0.2\textwidth]{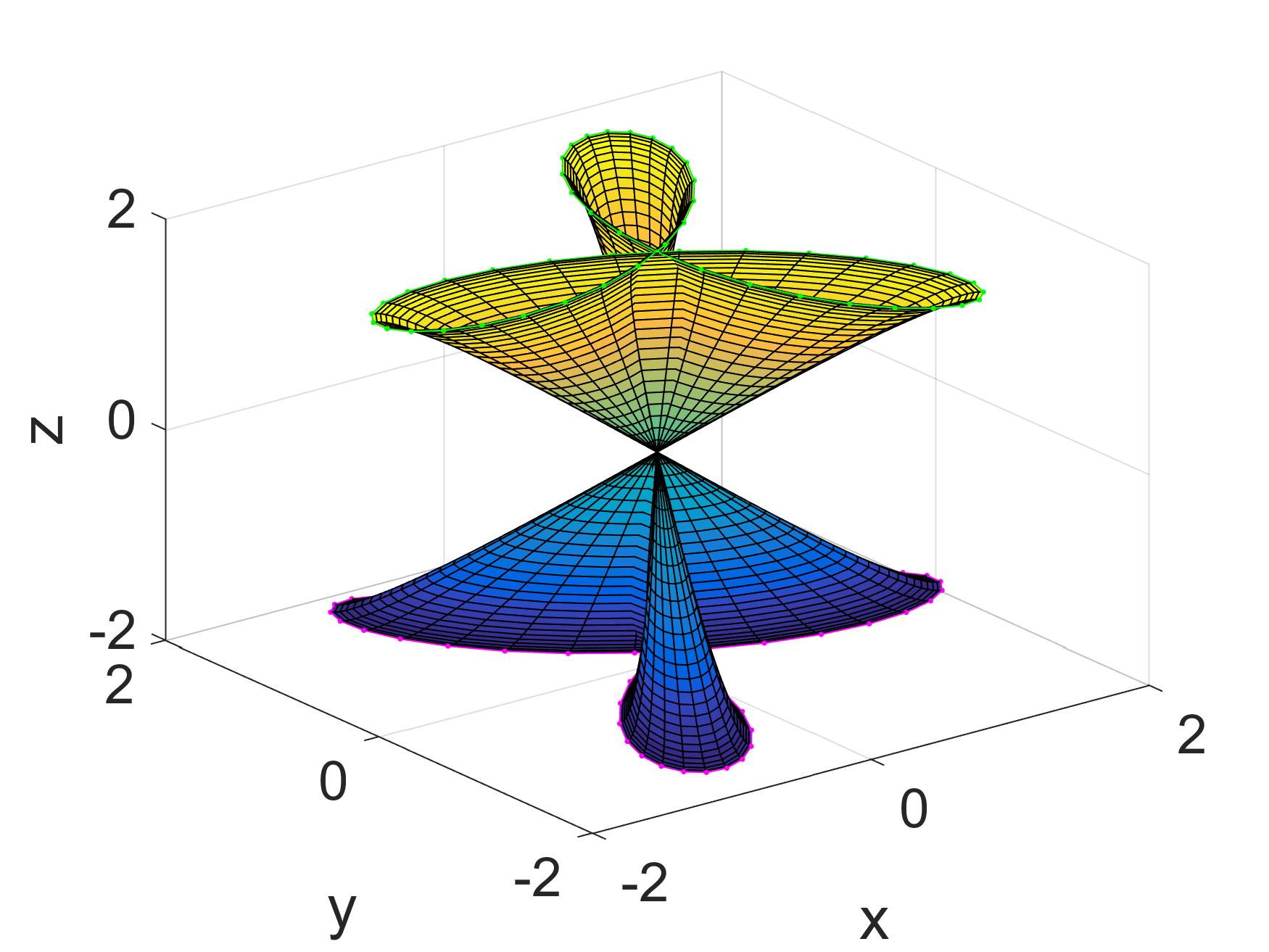}}
                \subfigure[]{\includegraphics[width=0.2\textwidth]{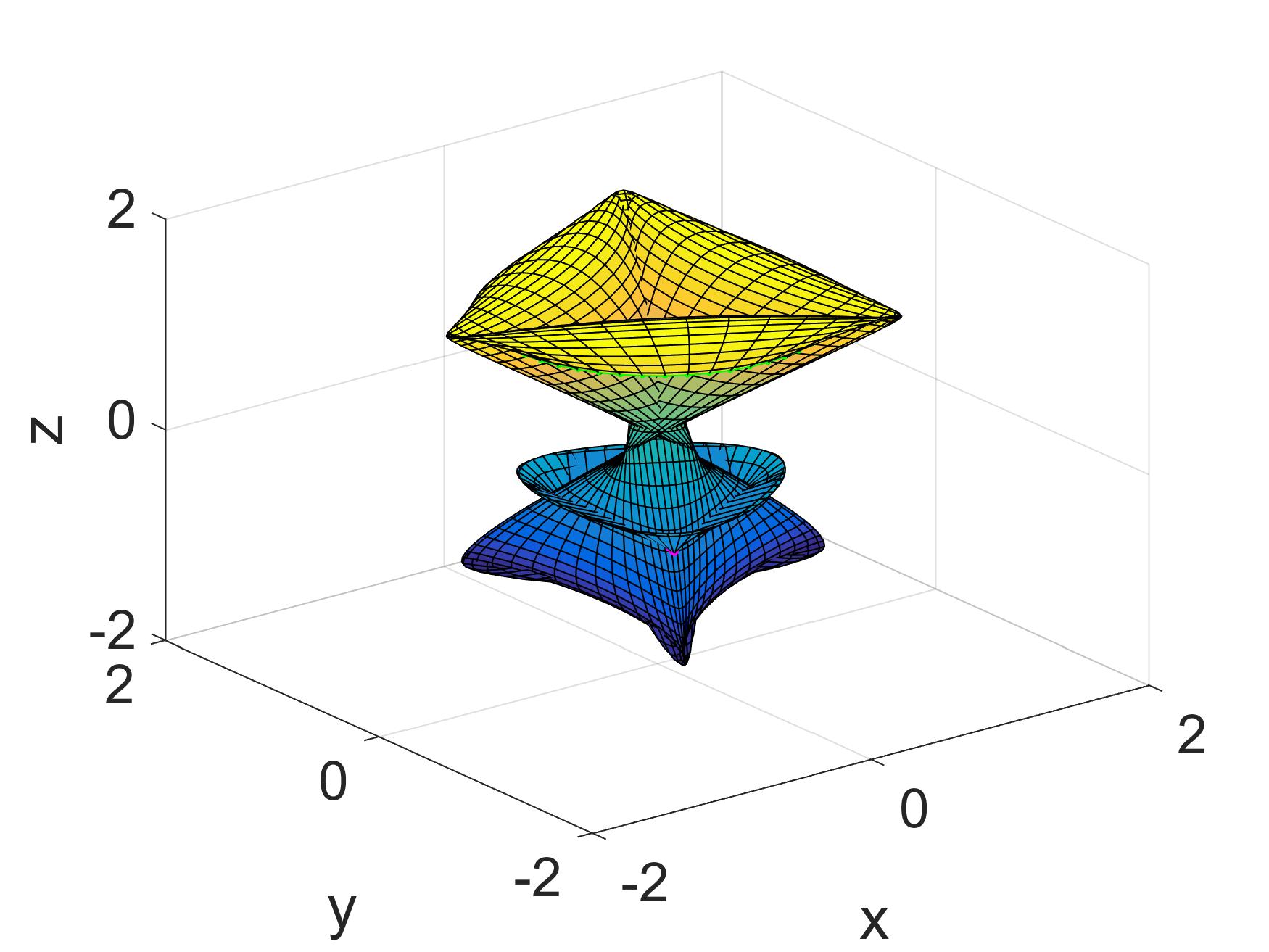}}\\
                \vspace{-0.41cm}  
                \subfigure[]{\includegraphics[width=0.2\textwidth]{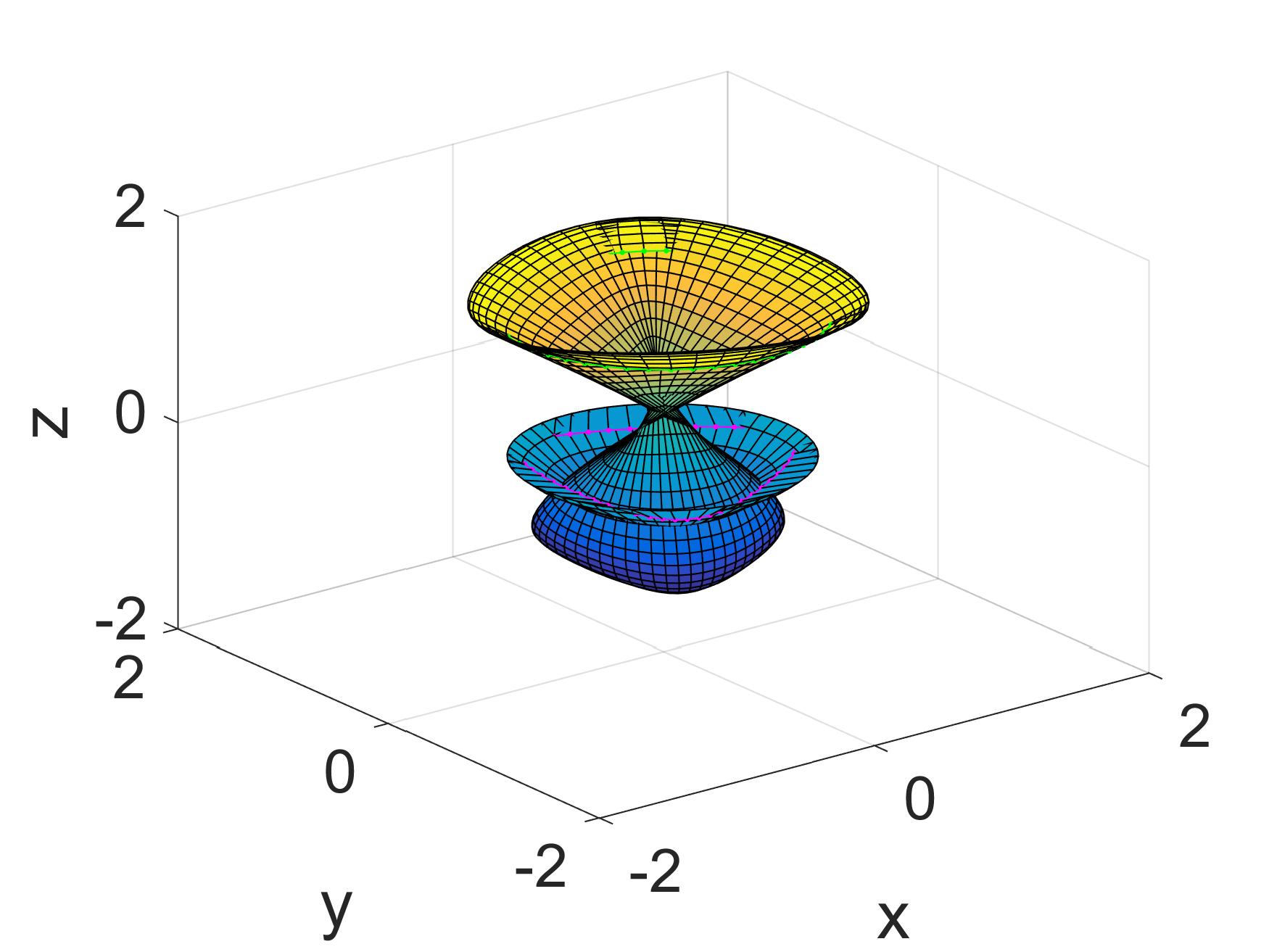}}
                \subfigure[]{\includegraphics[width=0.2\textwidth]{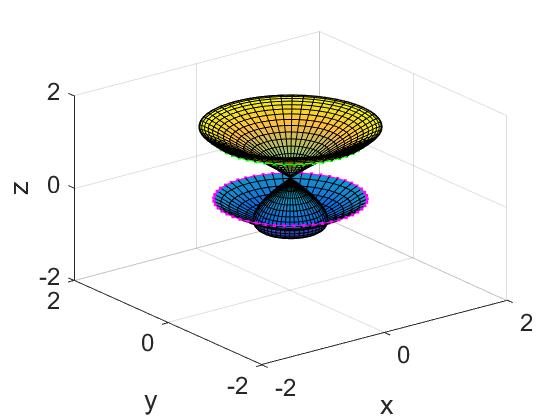}}\\
                \vspace{-0.41cm}
                \subfigure[]{\includegraphics[width=0.2\textwidth]{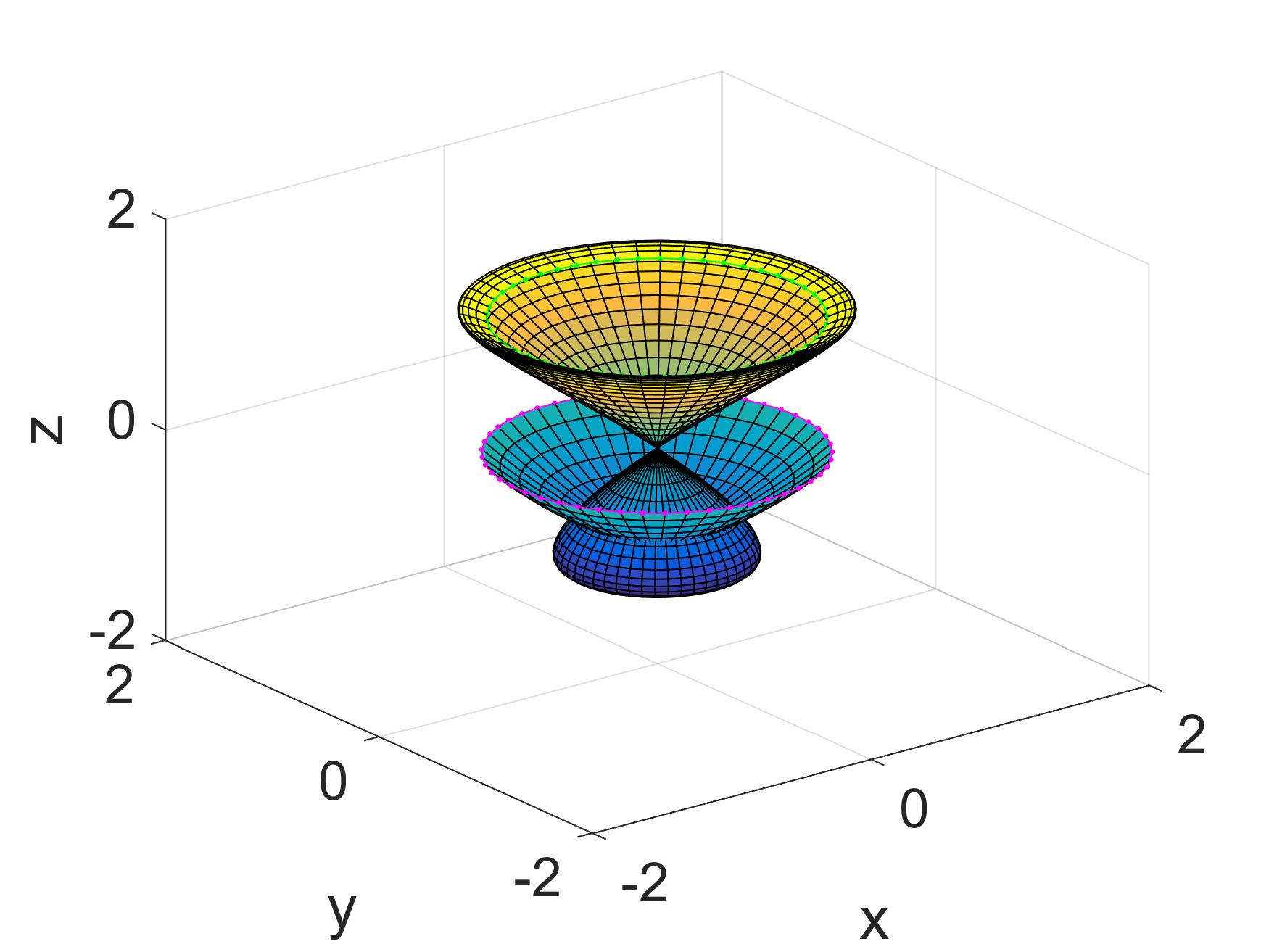}} 
                \subfigure[]{\includegraphics[width=0.2\textwidth]{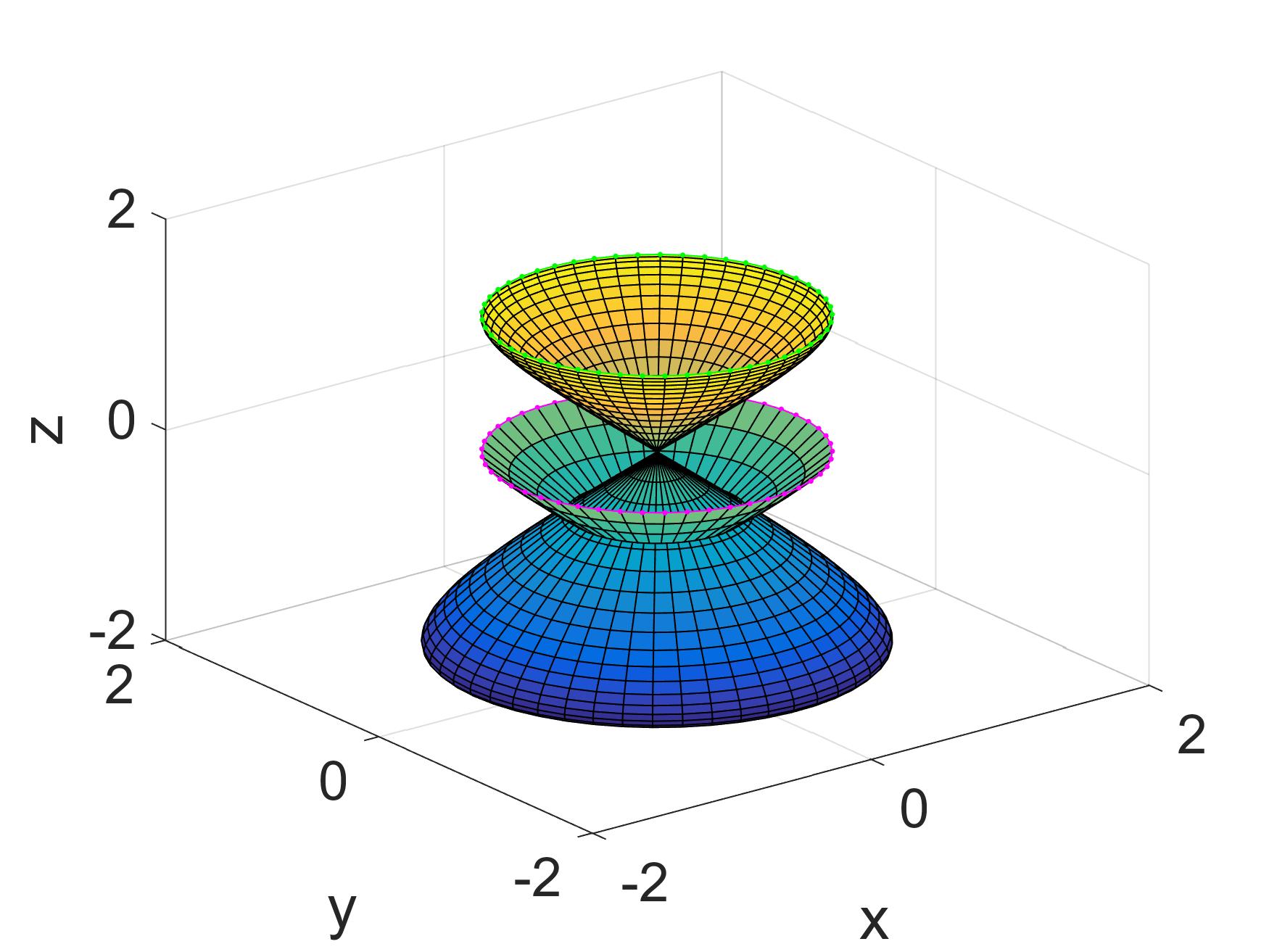}}
                \vspace{-0.21cm}
                \caption{The adaptive formation change process of the multi-agent system with unknown delay initial value $\hat D(0)=4$. (a) $t=0s$ (b) $t=0.09s$  (c) $t=0.2s$ \textcolor{black}{(d) $t=2s$ (e) $t=4s$} (f) $t=40s$.}\label{fig:6-formation}
        \end{figure}

\begin{figure}[htbp]
                \centering
                \subfigure[]{\includegraphics[width=0.32\textwidth]{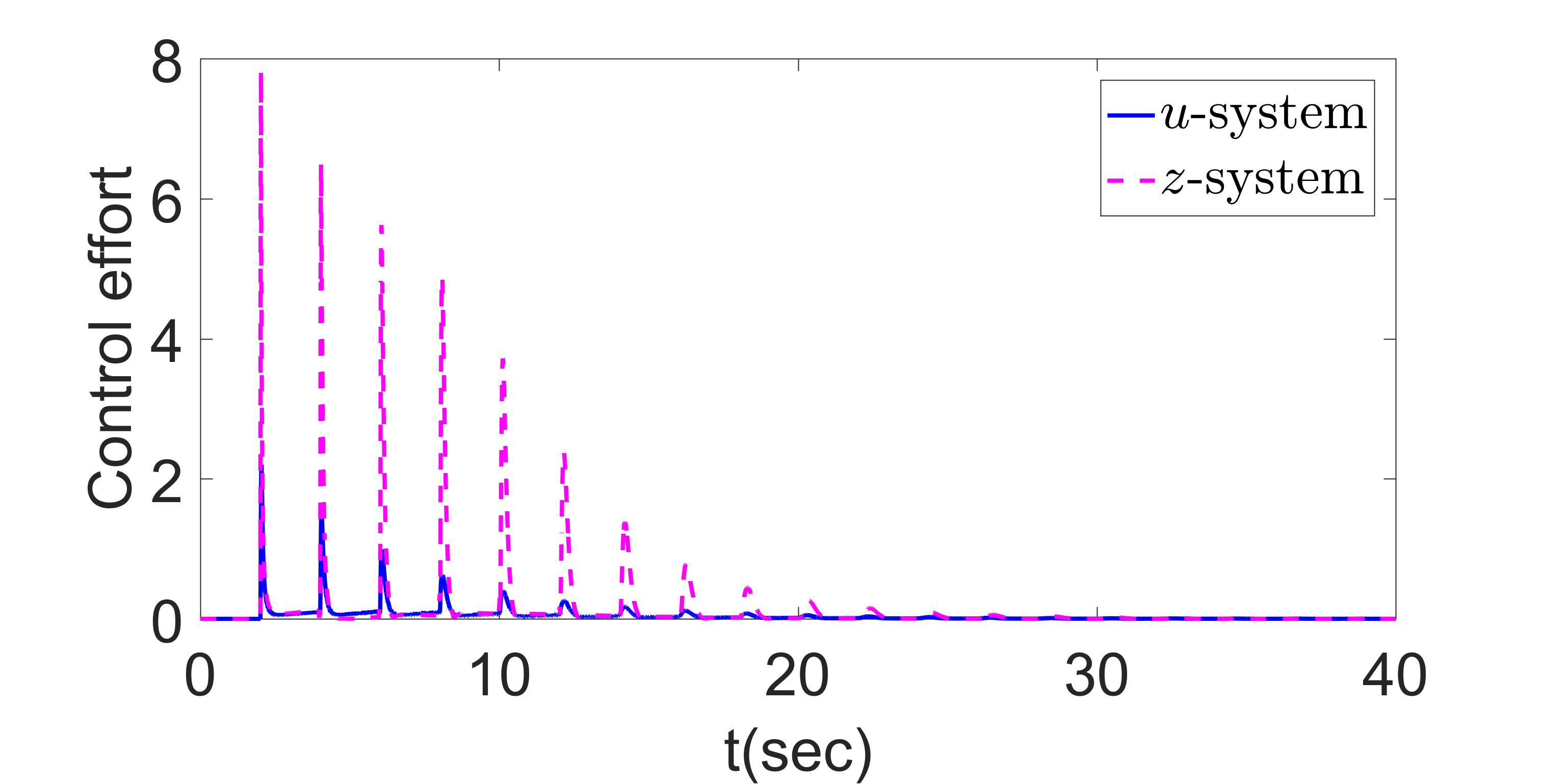}}
                \vspace{-0.21cm}
                \caption{Time-evolution of the control signals}\label{fig:6-control}
        \end{figure}

\begin{figure}[t]
                \centering 
                \subfigure[]{\includegraphics[width=0.23\textwidth]{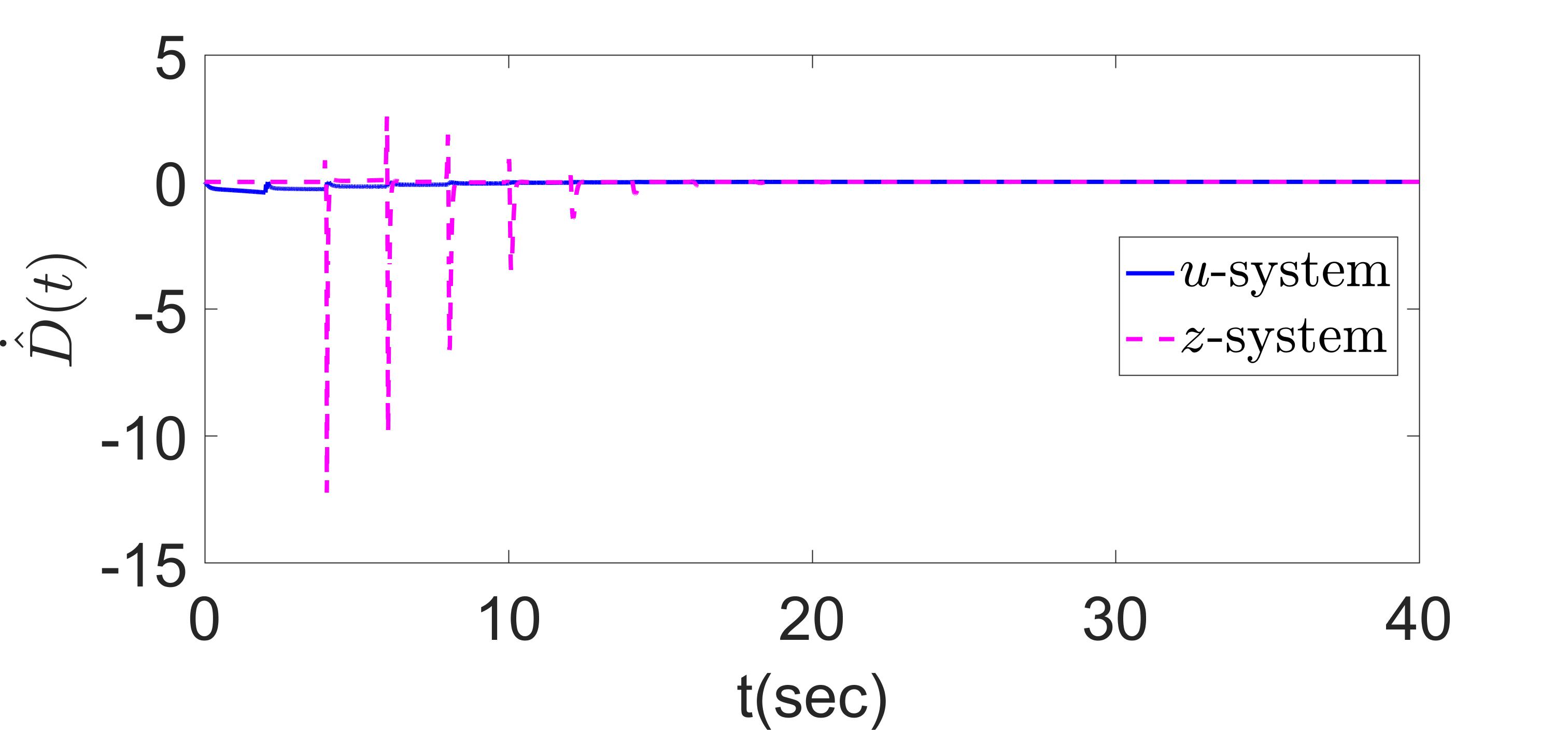}}
                \subfigure[]{\includegraphics[width=0.23\textwidth]{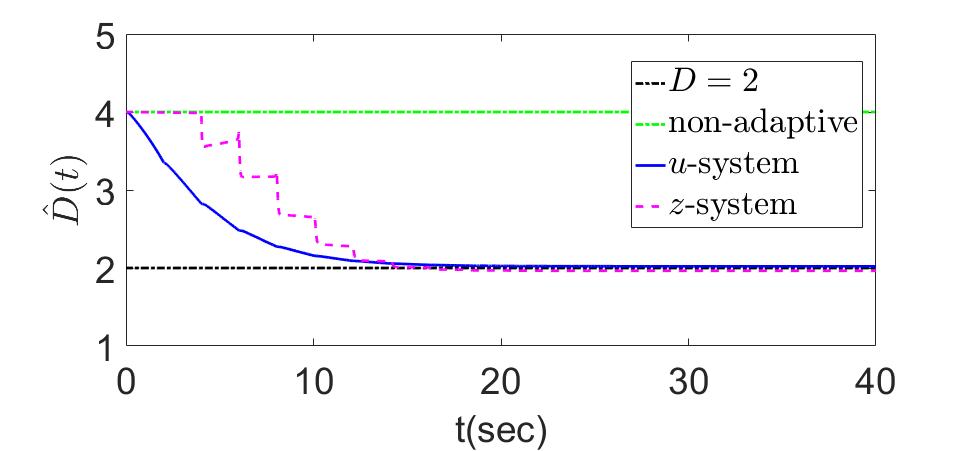}}
                \vspace{-0.21cm}
                \caption{  Delay estimate. (a) Dynamics of the updated law $\dot{\hat D}(t)$ (b) Time-evolution of the estimate of the unknown parameter $\hat D(t)$.}\label{figure2}
        \end{figure}

\begin{figure}[t]
                \centering
                \subfigure[]{\includegraphics[width=0.23\textwidth]{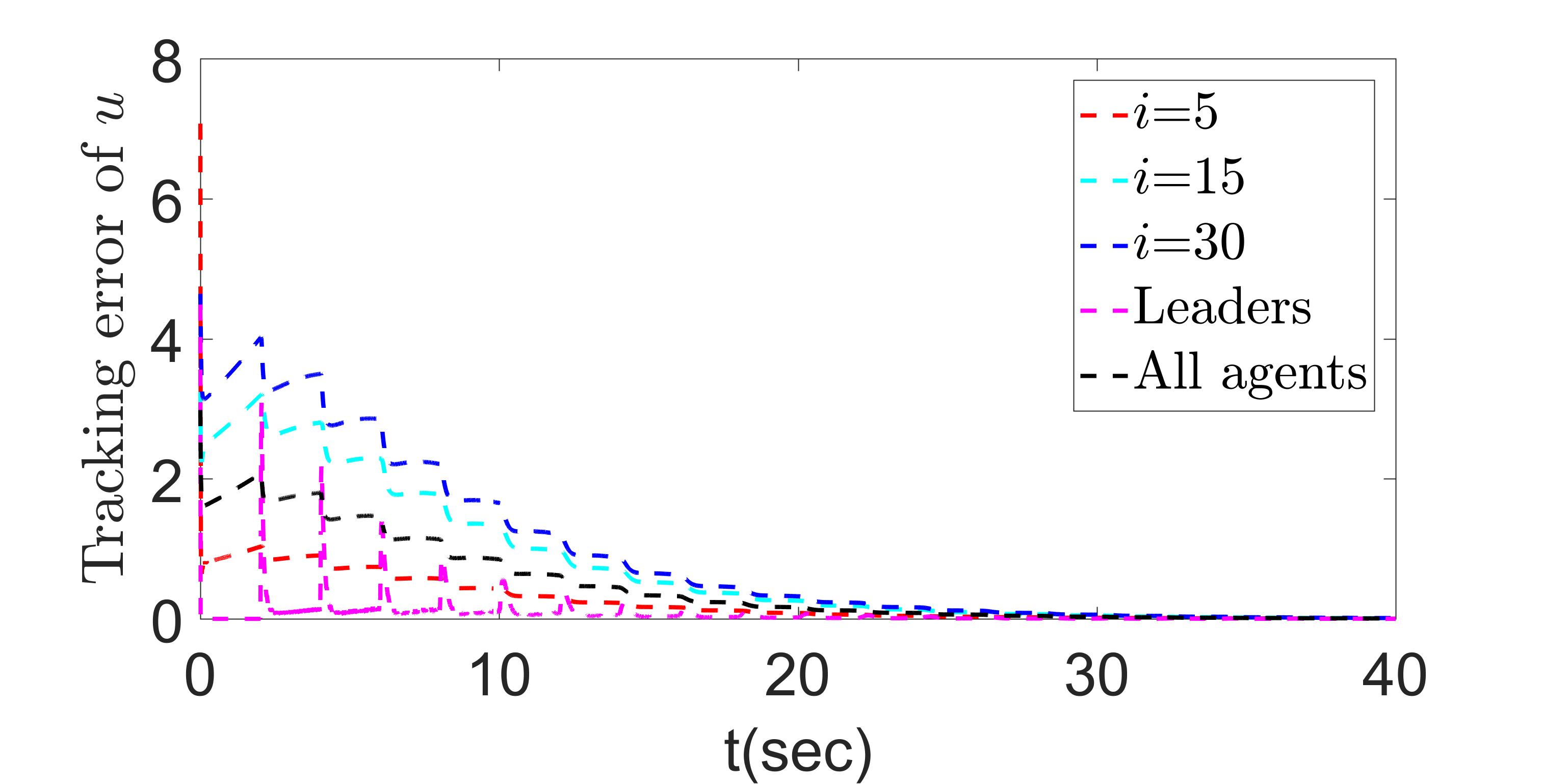}}
                \subfigure[]{\includegraphics[width=0.23\textwidth]{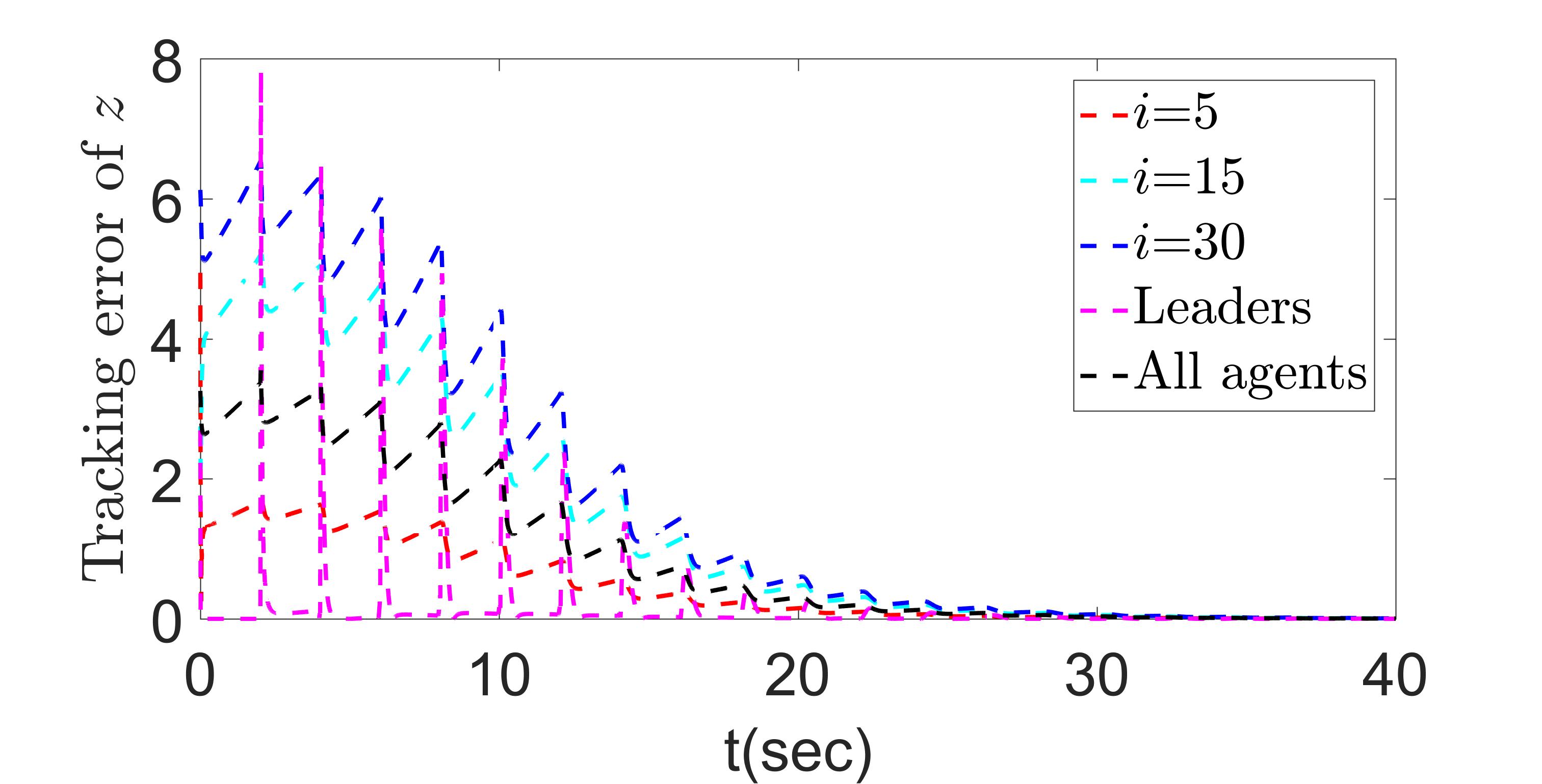}}\\
                \vspace{-0.41cm}
                \subfigure[]{\includegraphics[width=0.23\textwidth]{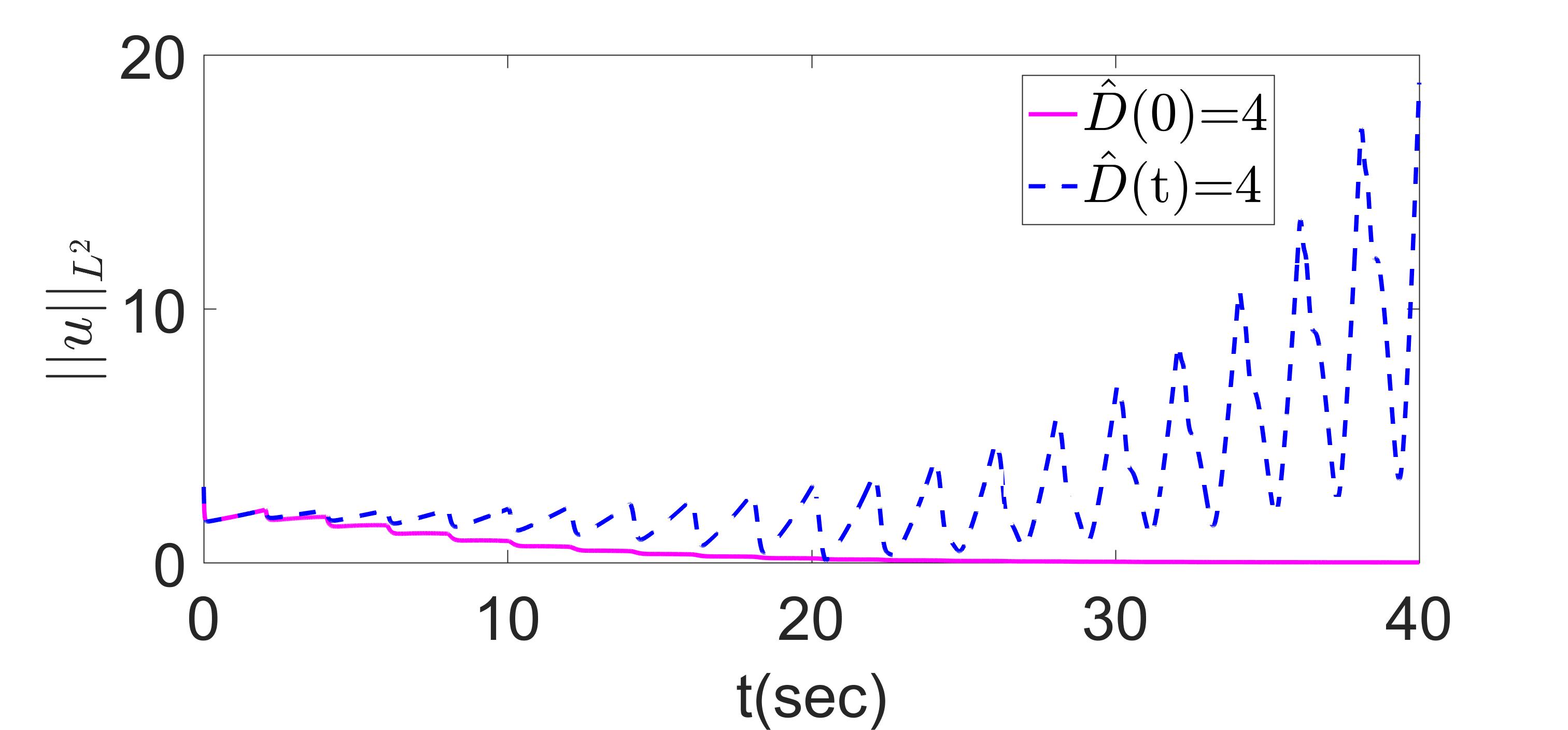}}
                \subfigure[]{\includegraphics[width=0.23\textwidth]{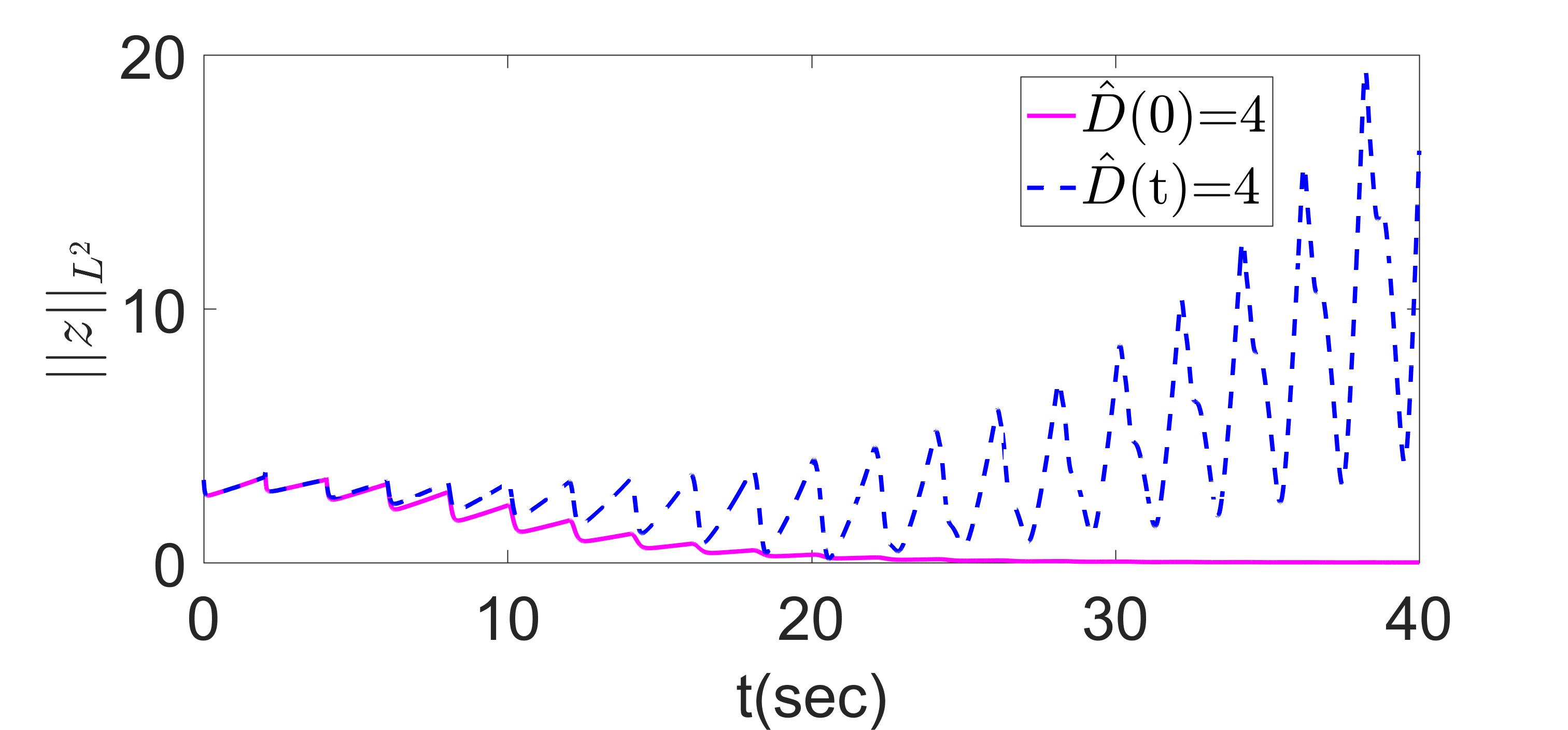}}
                \vspace{-0.21cm}
                \caption{(a) Tracking error of $u$ system under adaptive control, (b) Observation error of $z$ system under adaptive control,
(c) Average tracking error of $u$ system, (d) Average tracking error of $z$ system.}\label{fig:6-error}
        \end{figure}

        A formation control simulation example with $51\times50$ agents on a mesh grid in the 3-D space illustrates the performance of the proposed control laws with unknown input delay. The real value of input delay $D=2$, and the upper and lower bounds of the unknown delay are $\underline D=0.1$ and $\overline D=4$, respectively. The adaptive gain is fixed at $\varrho=0.05$. The model's parameters are    $\lambda_1=\lambda_2=10$, $\beta_1=\beta_2=0$. The control goal is to drive the formation  of the agents  from an initial equilibrium state  characterized by the  boundary values $f_1(\theta)=-e^{\mathrm j\theta}+e^{-\mathrm j2\theta}$, $g_1(\theta)=e^{\mathrm j\theta}-e^{-\mathrm j2\theta}$, $f_2(\theta)=-1.9$, $g_2(\theta)=1.9$ and the parameters  $\lambda_1=\lambda_2=10$, $\beta_1=\beta_2=0$ to a desired  formation  with boundary $f_1(\theta)=g_1(\theta)=e^{\mathrm j\theta}$, $f_2(\theta)=0$, $g_2(\theta)=1.3$ and the parameters of $\lambda_1=30$, $\lambda_2=20$, $\beta_1=\beta_2=1$. 
        Figure \ref{fig:6-formation} shows the formation diagram (or snapshots of the evolution in time) of a 3D multi-agent formation with an initial value of the unknown delay estimate $\hat D=4$ and from the initial  to the desired formation. The six snapshots of the formation's state illustrate the smooth evolution of collective dynamics between two different reference formations when the input delays are unknown.
        Figure \ref{fig:6-control} shows the time-evolution of the control signals, and it is clear that  the control effort  tends to zero and  ensures the stability of the closed-loop system dynamic. In Figure \ref{figure2}, (a) shows the dynamics of the update rate of the unknown parameter, $\dot {\hat {D}}$, when its initial value is $\hat D(0)=4$. It is clear that the updated rate gradually tends to zero over time; (b) describes the estimate of the unknown input delay for the system subject to  the designed  adaptive control law  for a given initial value $\hat D(0)=4$: the estimated  delay $\hat D$ gradually converges to the real value  $D=2$.   
        In Figure \ref{fig:6-error}, (a) and (b) show the tracking error of agents indexed by $i=5$, $i=15$,  $i=30$, and $i=51$ (actuator leaders) and the average of all  agents on the horizontal and vertical directions, respectively, under non-adaptive boundary control. It can be seen from the figure that the tracking error gradually tends to $0$ with time evolution. Figures (c) and (d) show  the $L^2$-norm of average tracking error of  all the agents in the horizontal and vertical directions, respectively. It can be seen that if the estimate of unknown delay $\hat D$  does not match the true value of the delay $D=2$, namely  if a delay mismatch occurs,  the tracking error diverges. 
        
        \section{Conclusion}\label{perspec}
        This paper  studies the formation control of MAS with unknown input delay in 3D space via cylindrical topology. To accomplish the desired 3D formation with stable transitions, we propose an adaptive controller with the backstepping method. The update law for estimating the unknown parameter is established using the Lyapunov method. As the dimensionality increases, the complexity of the problem grows significantly. To address this, we introduce a Fourier series to transform the PDE describing the two-dimensional cylindrical communication topology into the sum of infinite one-dimensional systems. Subsequently,  we prove the local stability of the closed-loop system and the regulation of the states to zero by a rather intricate Lyapunov function.
        In future work, we will extend our research to the systems subject to both unknown plant coefficients and input delays.

        \appendices 
        \section{The Proof of Proposition \ref{proposition4-1}}\label{apA}
        To prove the {norm equivalence}  between the state of the error system  \eqref{equ-phi}--\eqref{initial-vartheta} and that of the target system \eqref{equ-am0-adp}--\eqref{initial-amh-adp}, $S_i$, $i=1,2$ are constructed {using the bound of the integral of kernels, for example, let's consider the $L^2$ norm of $\vartheta(s,\theta,t)$. From equation \eqref{equ-trans4}}, we get the following estimate
        \begin{align}\label{equ-vs}
                &\int_0^1\int_\pi^\pi | {\vartheta}(s,\theta,t)|^2\mathrm d\theta\mathrm ds\nonumber\\
                \nonumber\leq&3\rVert  h\rVert^2+3\int_0^1\int_{-\pi}^\pi\bigg(\int_{0}^1\int_{-\pi}^\pi  \eta(s,\tau,\theta,\psi,\hat D)w(\tau,{\theta},t)\\
                \nonumber&\mathrm{d}\psi\mathrm{d}\tau\bigg)^2\mathrm d\theta\mathrm ds+3\hat D^2 \int_0^1\int_{-\pi}^\pi\bigg(\int_{0}^s\int_{-\pi}^\pi   q(s,\tau,\theta,\psi,\\
                &\hat D)h(\tau,{\theta},t)\mathrm{d}\psi\mathrm{d}\tau \bigg)^2\mathrm d\theta\mathrm ds.
        \end{align}
        From \eqref{equ-rn2} and \eqref{equ-qn2}, we have
        \begin{align}
                & \eta(s,\tau,\theta,\psi,\hat D)={Q(s,\theta-\psi,\hat D)} F(s,\tau),\\
                &q(s,\tau,\theta,\psi,\hat D)={Q(s-\tau,\theta-\psi,\hat D)}G(s,\tau),
        \end{align}
        where 
        \begin{align}
                F(s,\tau)&=2\sum_{i=1}^{\infty}e^{-\hat Di^2\pi^2s}\mathrm{sin}(i\pi\tau)\int_{0}^{1}\mathrm{sin}(i\pi\xi)l(1,\xi)\mathrm{d}\xi,\\ G(s,\tau)&=-2\sum_{i=1}^{\infty}e^{-\hat Di^2\pi^2(s-\tau)}i\pi(-1)^n\int_{0}^{1}\mathrm{sin}(i\pi\xi)l(1,\xi)\mathrm{d}\xi. \end{align}
        For the second term  on the left of the  inequality \eqref{equ-vs},
        using Fourier series, the following relations  hold \textcolor{black}{(see \eqref{fourier-Phi})
        \begin{align}
        w(s,\theta,t)&=\sum_{n=-\infty}^{\infty}w_n(s,t)e^{\mathrm{j}n\theta}\label{a6},
        \end{align}
        and \begin{align}
                {Q(s,\theta-\psi,\hat D)}&=\frac{1}{2\pi}\sum_{n=-\infty}^{\infty}e^{-\hat Dn^2s}e^{\mathrm{j}n(\theta-\psi)}.\label{a5}
        \end{align}}
        {It follows that
                \begin{align}
                        \nonumber&\int_{-\pi}^\pi  Q(s,\theta-\psi,\hat D)w(\tau,{\psi},t)\mathrm{d}\psi\\
                        \nonumber&=\frac{1}{2\pi}\sum_{n=-\infty}^{\infty}\sum_{m=-\infty}^{\infty}e^{-\hat Dn^2s}w_m(\tau,t)\int_{-\pi}^\pi e^{\mathrm{j}n(\theta-\psi)}e^{\mathrm{j}m\psi}\mathrm{d}\psi\\
                        &=\sum_{n=-\infty}^{\infty}e^{-\hat Dn^2s}e^{\mathrm{j}n\theta}w_n(\tau,t),
                \end{align}
                where we have used the orthogonality property of the Fourier series (the same conclusion is reached using the convolution theorem).
        }
        From the  Parseval's theorem, the following can be deduced
        \begin{align}
                \nonumber&\int_{-\pi}^\pi\bigg( \int_{-\pi}^\pi  Q(s,\theta-\psi,\hat D)w(\tau,{\psi},t)\mathrm{d}\psi F(s,\tau)\bigg)^2\mathrm d\theta\\
                \nonumber&=2\pi\sum_{n=-\infty}^{\infty}e^{-2\hat Dn^2s}|w_n(\tau,t)|^{2}|F(s,\tau)|^{2}\\
                &{\leq 2\pi\sum_{n=-\infty}^{\infty}|w_n(\tau,t)|^{2}|F(s,\tau)|^{2},}
        \end{align}
        which allow to state that
        \begin{align}
                \nonumber&\int_0^1\int_{-\pi}^\pi\bigg( \int_0^1\int_{-\pi}^\pi  Q(s,\theta-\psi,\hat D)F(s,\tau)\\
                \nonumber&\cdot w(\tau,{\psi},t)\mathrm{d}\psi\mathrm d\tau\bigg)^2\mathrm d\theta\mathrm ds\\
                \nonumber&\leq\int_0^1\bigg(\int_{-\pi}^\pi \int_0^1\bigg(\int_{-\pi}^\pi  Q(s,\theta-\psi,\hat D)w(\tau,{\psi},t)\mathrm{d}\psi\bigg)^2\mathrm d\tau\mathrm d\theta\\
                \nonumber&\cdot \int_0^1|F(s,\tau)|^2\mathrm d\tau\bigg)\mathrm ds\\
                \nonumber&\leq\int_0^1 \int_0^1\int_{-\pi}^\pi\bigg(\int_{-\pi}^\pi  Q(s,\theta-\psi,\hat D)w(\tau,{\psi},t)\mathrm{d}\psi\bigg)^2\mathrm d\tau\mathrm d\theta\mathrm ds\\
                \nonumber&\cdot \int_0^1\int_0^1|F(s,\tau)|^2\mathrm d\tau\mathrm ds\\
                \nonumber&\leq{2\pi}\int_0^1\sum_{n=-\infty}^{\infty}|w_n(\tau,t)|^{2}\mathrm d\tau\int_0^1\int_0^1|F(s,\tau)|^2\mathrm d\tau\mathrm ds\\
                &=\int_0^1\int_0^1|F(s,\tau)|^2\mathrm d\tau\mathrm ds\rVert w\rVert^2.\label{equ-second-term}
        \end{align}
        Similarly, for the third term \textcolor{black}{on the left of the  inequality \eqref{equ-vs}, based on  \eqref{a5}, \eqref{a6} and the orthogonality property of Fourier series, one can get}
        \begin{align}
                \nonumber&\int_{-\pi}^\pi  Q(s-\tau,\theta-\psi,\hat D)h(\tau,{\psi},t)\mathrm{d}\psi\\
                &=\sum_{n=-\infty}^{\infty}e^{-\hat Dn^2(s-\tau)}e^{\mathrm{j}n\theta}h_n(\tau,t).
        \end{align}
        From the Parseval's theorem, the following can be deduced
        \begin{align}
                \nonumber&\int_{-\pi}^\pi\bigg( \int_{-\pi}^\pi  Q(s-\tau,\theta-\psi,\hat D)G(s,\tau)h(\tau,{\psi},t)\mathrm{d}\psi \bigg)^2\mathrm d\theta\\
                \nonumber&=2\pi\sum_{n=-\infty}^{\infty}e^{-2\hat Dn^2(s-\tau)}|h_n(\tau,t)|^{2}|G(s,\tau)|^{2}\\
                &{\leq 2\pi\sum_{n=-\infty}^{\infty}|h_n(\tau,t)|^{2}|G(s,\tau)|^{2},}
        \end{align}
        and consequently
        \begin{align}
                \nonumber&\int_0^1\int_{-\pi}^\pi\bigg( \int_0^s\int_{-\pi}^\pi  Q(s-\tau,\theta-\psi,\hat D)G(s,\tau)\\
                \nonumber&\cdot h(\tau,{\psi},t)\mathrm{d}\psi\mathrm d\tau\bigg)^2\mathrm d\theta\mathrm ds\\
                \nonumber&\leq\int_0^1 \int_0^s\int_{-\pi}^\pi\bigg(\int_{-\pi}^\pi  Q(s-\tau,\theta-\psi,\hat D)\\
                \nonumber&\cdot h(\tau,{\psi},t)\mathrm{d}\psi\bigg)^2\mathrm d\tau\mathrm d\theta\mathrm ds\int_0^1\int_0^s|G(s,\tau)|^2\mathrm d\tau\mathrm ds\\
                \nonumber&\leq{2\pi}\int_0^1\sum_{n=-\infty}^{\infty}|h_n(\tau,t)|^{2}\mathrm d\tau\int_0^1\int_0^s|G(s,\tau)|^2\mathrm d\tau\mathrm ds\\
                &\leq \int_0^1\int_0^s|G(s,\tau)|^2\mathrm d\tau\mathrm ds\rVert h\rVert^2.\label{equ-third-term}
        \end{align}
        Thus, \textcolor{black}{combining with \eqref{equ-second-term} and \eqref{equ-third-term},} we get
        \begin{align}\label{equ-vs0}
                &\int_0^1\int_\pi^\pi | {\vartheta}(s,\theta,t)|^2\mathrm d\theta\mathrm ds\leq3\rVert  h\rVert^2\nonumber\\
                \nonumber&+3\int_0^1\int_{-\pi}^\pi\bigg(\int_{0}^1\int_{-\pi}^\pi  {Q(s,\theta-\psi,\hat D)} F(s,\tau)w(\tau,{\theta},t)\\
                \nonumber&\mathrm{d}\psi\mathrm{d}\tau\bigg)^2\mathrm d\theta\mathrm ds+3\hat D^2 \int_0^1\int_{-\pi}^\pi\bigg(\int_{0}^s\int_{-\pi}^\pi   {Q(s-\tau,\theta-\psi,\hat D)}\\
                \nonumber&\cdot G(s,\tau)h(\tau,{\theta},t)\mathrm{d}\psi\mathrm{d}\tau \bigg)^2\mathrm d\theta\mathrm ds\\
                \nonumber\leq&3(1+\overline D^2)\int_0^1\int_0^s|G(s,\tau)|^2\mathrm d\tau\mathrm ds\rVert  h\rVert^2\\
                &+3\int_0^1\int_0^1|F(s,\tau)|^2\mathrm d\tau\mathrm ds\rVert w\rVert^2,
        \end{align}
        where $\int_0^1\int_0^s|G(s,\tau)|^2\mathrm d\tau\mathrm ds$ and $\int_0^1\int_0^1|F(s,\tau)|^2\mathrm d\tau\mathrm ds$ are  bounded as established in \cite{WANG2021109909}, which complete the proof of Proposition \ref{proposition4-1}. 
        
        \bibliographystyle{plain}       
        \bibliography{multi_agent_delay_ref}

\begin{thebibliography}{10}

\bibitem{ai2021distributed}
X.~Ai and L.~Wang.
\newblock Distributed fixed-time event-triggered consensus of linear
  multi-agent systems with input delay.
\newblock {\em International Journal of Robust and Nonlinear Control},
  31(7):2526--2545, 2021.

\bibitem{alonso2019distributed}
J.~Alonso-Mora, E.~Montijano, T.~N{\"a}geli, O.~Hilliges, M.~Schwager, and
  D.~Rus.
\newblock Distributed multi-robot formation control in dynamic environments.
\newblock {\em Autonomous Robots}, 43(5):1079--1100, 2019.

\bibitem{Mora2015}
J.~Alonso-Mora, T.~Naegeli, P.~Beardsley, and P.~Beardsley.
\newblock Collision avoidance for aerial vehicles in multi-agent scenarios.
\newblock {\em Autonomous Robots}, 39(1):101--121, 2015.

\bibitem{Brown2009}
J.~W. Brown and R.~V. Churchill.
\newblock {\em Complex variables and applications}.
\newblock Brown and Churchill series. McGraw-Hill Higher Education, 2009.

\bibitem{Fax2004}
J.~A. Fax and R.~M. Murray.
\newblock Information flow and cooperative control of vehicle formations.
\newblock {\em IEEE Transactions on Automatic Control}, 49(9):1465--1476, 2004.

\bibitem{Ferrari2006}
G.~Ferrari-Trecate, A.~Buffa, and M.~Gati.
\newblock Analysis of coordination in multi-agent systems through partial
  difference equations.
\newblock {\em IEEE Transactions on Automatic Control}, 51(6):1058--1063, 2006.

\bibitem{FREUDENTHALER2020108897}
G.~Freudenthaler and T.~Meurer.
\newblock {PDE}-based multi-agent formation control using flatness and
  backstepping: Analysis, design and robot experiments.
\newblock {\em Automatica}, 115:108897, 2020.

\bibitem{Frihauf2011}
P.~Frihauf and M.~Krstic.
\newblock Leader-enabled deployment onto planar curves: A {P}{D}{E}-based
  approach.
\newblock {\em IEEE Transactions on Automatic Control}, 56(8):1791--1806, 2011.

\bibitem{hou2017consensus}
W.~Hou, M.~Fu, H.~Zhang, and Z.~Wu.
\newblock Consensus conditions for general second-order multi-agent systems
  with communication delay.
\newblock {\em Automatica}, 75:293--298, 2017.

\bibitem{karafyllis2019adaptive}
I.~Karafyllis, M.~Kontorinaki, and M.~Krstic.
\newblock Adaptive control by regulation-triggered batch least squares.
\newblock {\em IEEE Transactions on Automatic Control}, 65(7):2842--2855, 2019.

\bibitem{MKrstic2009}
M.~{Krstic}.
\newblock Control of an unstable reaction-diffusion {PDE} with long input
  delay.
\newblock {\em Systems \& Control Letters}, 58(10):773--782, 2009.

\bibitem{Lee2006}
D.~Lee and M.~W. Spong.
\newblock Agreement with non-uniform information delays.
\newblock In {\em American Control Conference (ACC)}, pages 756--761, 2006.

\bibitem{9107091}
K.~Li, C.~Hua, X.~You, and X.~Guan.
\newblock Distributed output-feedback consensus control for nonlinear
  multiagent systems subject to unknown input delays.
\newblock {\em IEEE Transactions on Cybernetics}, 52(2):1292--1301, 2022.

\bibitem{Lin2014}
P.~Lin and W.~Ren.
\newblock Constrained consensus in unbalanced networks with communication
  delays.
\newblock {\em IEEE Transactions on Automatic Control}, 59(3):775--781, 2014.

\bibitem{LIU2021104927}
Z.~Liu, D.~Nojavanzadeh, D.~Saberi, A.~Saberi, and A.~A. Stoorvogel.
\newblock Scale-free protocol design for regulated state synchronization of
  homogeneous multi-agent systems with unknown and non-uniform input delays.
\newblock {\em Systems \& Control Letters}, 152:104927, 2021.

\bibitem{Meurer2011}
T.~Meurer and M.~Krstic.
\newblock Finite-time multi-agent deployment: A nonlinear {P}{D}{E} motion
  planning approach.
\newblock {\em Automatica}, 47(11):2534--2542, 2011.

\bibitem{Qi2015}
J.~Qi, R.~Vazquez, and M.~Krstic.
\newblock Multi-agent deployment in 3-{D} via {P}{D}{E} control.
\newblock {\em IEEE Transactions on Automatic Control}, 60(4):891--906, 2015.

\bibitem{qi2019control}
J.~Qi, S.~Wang, J.~Fang, and M.~Diagne.
\newblock Control of multi-agent systems with input delay via {PDE}-based
  method.
\newblock {\em Automatica}, 106:91--100, 2019.

\bibitem{qi2017wave}
J.~Qi, J.~Zhang, and Y.~Ding.
\newblock Wave equation-based time-varying formation control of multiagent
  systems.
\newblock {\em IEEE Transactions on Control Systems Technology},
  26(5):1578--1591, 2018.

\bibitem{Krstic2010}
A.~Smyshlyaev and M.~Krstic.
\newblock {\em Adaptive Control of Parabolic {PDE}s}.
\newblock Princeton University Press, 2010.

\bibitem{tan2017leader}
X.~Tan, J.~Cao, X.~Li, and A.~Alsaedi.
\newblock Leader-following mean square consensus of stochastic multi-agent
  systems with input delay via event-triggered control.
\newblock {\em IET Control Theory \& Applications}, 12(2):299--309, 2017.

\bibitem{tang2017formation}
S.~Tang, J.~Qi, and J.~Zhang.
\newblock Formation tracking control for multi-agent systems: A wave-equation
  based approach.
\newblock {\em International Journal of Control, Automation and Systems},
  15(6):2704--2713, 2017.

\bibitem{Tian2008}
Y.~Tian and C.~Liu.
\newblock Consensus of multi-agent systems with diverse input and communication
  delays.
\newblock {\em IEEE Transactions on Automatic Control}, 53(9):2122--2128, 2008.

\bibitem{vazquez2016explicit}
R.~Vazquez and M.~Krstic.
\newblock Explicit output-feedback boundary control of reaction-diffusion
  {P}{D}{E}s on arbitrary-dimensional balls.
\newblock {\em ESAIM: Control, Optimisation and Calculus of Variations},
  22(4):1078--1096, 2016.

\bibitem{Wang2017MAS}
H.~Wang, D.~Guo, X.~Liang, W.~Chen, G.~Hu, and K.~K. Leang.
\newblock Adaptive vision-based leader-follower formation control of mobile
  robots.
\newblock {\em IEEE Transactions on Industrial Electronics}, 64(4):2893--2902,
  2017.

\bibitem{Jiwang2023delay}
J.~Wang and M.~Diagne.
\newblock Delay-adaptive boundary control of coupled hyperbolic {PDE}-{ODE}
  cascade systems.
\newblock {\em arXiv e-prints}, pages arXiv--2301, 2023.

\bibitem{wang2021delay}
S.~Wang, M.~Diagne, and J.~Qi.
\newblock Delay-adaptive predictor feedback control of
  reaction--advection--diffusion {PDEs} with a delayed distributed input.
\newblock {\em IEEE Transactions on Automatic Control}, 67(7):3762--3769, 2022.

\bibitem{WANG2021109909}
S.~Wang, J.~Qi, and M.~Diagne.
\newblock Adaptive boundary control of reaction–diffusion {PDE}s with unknown
  input delay.
\newblock {\em Automatica}, 134:109909, 2021.

\bibitem{Wang2017}
S.~{Wang}, J.~{Qi}, and J.~{Fang}.
\newblock Control of 2-{D} reaction-advection-diffusion {PDE} with input delay.
\newblock In {\em 2017 Chinese Automation Congress (CAC)}, pages 7145--7150,
  2017.

\bibitem{wang2023delay}
S.~Wang, J.~Qi, and M.~Krstic.
\newblock Delay-adaptive control of first-order hyperbolic {PIDEs}.
\newblock {\em arXiv preprint arXiv:2307.04212}, 2023.

\bibitem{yu2010some}
W.~Yu, G.~Chen, and M.~Cao.
\newblock Some necessary and sufficient conditions for second-order consensus
  in multi-agent dynamical systems.
\newblock {\em Automatica}, 46(6):1089--1095, 2010.

\bibitem{Zetocha2000}
P.~Zetocha, L.~Self, R.~Wainwright, R.~Burns, M.~Brito, and D.~Surka.
\newblock Commanding and controlling satellite clusters.
\newblock {\em IEEE Intelligent Systems and their Applications}, 15(6):8--13,
  2000.

\bibitem{Zhangunknowndelay2021}
M.~Zhang, A.~Saberi, and A.~A. Stoorvogel.
\newblock Semi-global state synchronization for multi-agent systems subject to
  actuator saturation and unknown nonuniform input delay.
\newblock {\em IEEE Transactions on Network Science and Engineering},
  8(1):488--497, 2021.

\bibitem{zhu2015event}
W.~Zhu and Z.~Jiang.
\newblock Event-based leader-following consensus of multi-agent systems with
  input time delay.
\newblock {\em IEEE Transactions on Automatic Control}, 60(5):1362--1367, 2015.

\end{thebibliography}

        
        
        
\end{document}